\documentclass[10pt] {amsart}
\usepackage{amsfonts}
\usepackage[] {amsmath}
\usepackage[] {amssymb}
\usepackage[] {amsthm}
\usepackage[] {amsfonts}

\newtheorem {theorem}    {Theorem}    [section]
\newtheorem {lemma}      {Lemma}      [section]
\newtheorem {definition} {Definition}   [section]
\newtheorem {corollary}  {Corollary}    [section]

\newtheorem* {conjecture*}       {Conjecture}
\newtheorem* {theorem*}          {Theorem}

\newtheorem* {acknowledgements*} {Acknowledgements}

\theoremstyle {definition} \newtheorem {remark}  {Remark}  [section]
\theoremstyle {definition} \newtheorem {example} {Example} [section]

\newcommand \Tan    {\operatorname{Tan}}

\newcommand \R      {\operatorname{R}}
\newcommand \F       {\operatorname {F}}
\newcommand \inj    {\operatorname{inj}}
\newcommand \sing   {\operatorname{sing}}
\newcommand \reg    {\operatorname{reg}}
\newcommand \spt    {\operatorname{spt}}

\newcommand \tr     {\operatorname{tr}}
\renewcommand \div  {\operatorname{div}}
\newcommand \dist   {\operatorname{dist}}
\newcommand \diam   {\operatorname{diam}}
\newcommand \mc     {\operatorname{H}}
\newcommand \D      {\operatorname{D}}
\renewcommand \L      {\operatorname{L}}
\newcommand \Rad      {\operatorname{Rad}}
\newcommand \B {\operatorname {B}}

\newcommand {\Ric}  {\operatorname{Ric}}
\newcommand {\graph}{\operatorname{graph}}

\title [The Plateau problem] {The Plateau problem for marginally outer trapped surfaces}
\author{Michael Eichmair}

\address {Michael Eichmair, Department of Mathematics, Massachusetts Institute of Technology, Cambridge MA-02139, USA}
\email{eichmair@math.mit.edu}

\begin{document}

\maketitle

\begin{abstract} We solve the Plateau problem for marginally outer trapped surfaces in general Cauchy data sets. We employ the Perron method and tools from geometric measure theory to force and control a blow-up of Jang's equation. Substantial new geometric insights regarding the lower order properties of marginally outer trapped surfaces are gained in the process. The techniques developed in this paper are flexible and can be adapted to other non-variational existence problems.
\end{abstract}

\section{Introduction}  \label{sec: introduction}

Let $(M^n, g)$ be a complete oriented Riemannian manifold and let $p$ be a
symmetric $(0, 2)$-tensor field defined on $M^n$. In the context of general
relativity, triples $(M^n, g, p)$ of this form arise as spacelike hypersurfaces
of time-oriented Lorentz manifolds $(\bar{M}^{n + 1}, \bar{g})$ where $g$ is
the induced metric and $p$ is the second fundamental form with respect to the
future pointing normal of $M^n$ in $\bar{M}^{n+1}$. We will refer to $(M^n, g, p)$ as an initial data set of general relativity. The case classically studied
by physicists is that of $n=3$, the development of string theory however has greatly increased the interest in higher dimensional space-times. Since the 
methods presented in this paper are in no way limited to $3$ dimensions and, we hope, interesting in their own right, we will keep the exposition general. The author is grateful to G. Galloway for suggesting that the case $n=2$ should also be included in the discussion because of interest in quantum gravity.   \\

We now fix the terminology used in this work and briefly indicate the physical significance of the surfaces whose existence and properties are studied in this paper. We refer to \cite{HawEll}, \cite{Wald}, the discussion in \cite{AM}, and the survey article \cite{ChruscielPollackGalloway} for further background and concise references. Given an initial data set $(M^n, g, p)$ and an open subset $\Omega \subset M^n$ we say that a smooth embedded oriented hypersurface $\Sigma^{n-1} \subset \Omega$, possibly with boundary, is \emph{marginally outer trapped} if it satisfies the prescribed mean curvature equation \begin{equation} \label{eq: marginally trapped surface equation} \mc_{\Sigma} + \tr_{\Sigma}(p) = 0. \end{equation}
Here, $\tr_{\Sigma} (p)$ denotes the trace of the tensor $p$ over the tangent space of $\Sigma^{n-1}$ and $\mc_\Sigma$ is the scalar mean curvature of $\Sigma^{n-1}$ computed consistently with the orientation of $\Sigma^{n-1}$. If $\Sigma^{n-1}$ is given specifically as a (relative) boundary in $\Omega$ we require that condition (\ref{eq: marginally trapped surface equation}) holds for the mean curvature scalar $\mc_\Sigma$ computed as the tangential divergence of the outward pointing unit normal of $\Sigma^{n-1}$. The conditions $\mc_\Sigma + \tr_\Sigma(p) < 0$ and $\mc_\Sigma + \tr_\Sigma (p) > 0$ for such a relative boundary $\Sigma^{n-1}$ are called \emph{outer  trapped} and \emph{outer untrapped} respectively.
If a space-time $(\bar{M}^{n+1}, \bar{g})$ containing $(M^n, g, p)$ as a Cauchy surface satisfies the null energy condition and is asymptotically predictable, then the singularity theorem of Hawking and Penrose implies that closed
marginally outer trapped surfaces $\Sigma^{n-1} \subset M^n$ lie inside the black hole region of $(\bar{M}^{n+1}, \bar{g})$. 
In this concrete sense, marginally outer trapped surfaces have been studied as a ``quasi-local version of the event horizon" of $\bar{M}^{n+1}$ traced out on $M^n$.
If $\Sigma^{n-1}$ is outermost with respect to inclusion in its homology class amongst all marginally outer trapped surfaces of $(M^n, g, p)$, then it is referred to as an \emph{apparent horizon} in the literature, and the compact region bounded by  $\Sigma^{n-1}$ is called the \emph{trapped region}. 
\\

The prescribed mean curvature condition $\mc_\Sigma + \tr_\Sigma (p) = 0$ characterizing marginally outer trapped surfaces is not known to arise as the
Euler-Lagrange equation of a variational problem on $(M^n, g, p)$, except for the cases where the data is totally geodesic (i.e. $p \equiv 0$) and $\Sigma^{n-1}$ is a minimal surface, or, more generally, when $p$ is a
constant multiple of $g$ and marginally outer trapped surfaces are those with constant mean curvature. \\

In their proof of the space-time positive mass theorem \cite{PMTII}, R. Schoen and S.-T. Yau constructed solutions to a one parameter family of regularizations of Jang's equation on an initial data set $(M^n, g, p)$. Using curvature estimates and the Harnack principle on the graphs of these solutions in $M^n \times \mathbb{R}$, they argued that they can be passed subsequentially in the regularization limit to complete embedded submanifolds of $M^n \times \mathbb{R}$ that are made up of graphical components, corresponding to actual solutions of Jang's equation on open subsets of $M^n$, and cylindrical components, whose cross-sections are closed marginally outer trapped surfaces in $M^n$. It is useful to think of this limiting submanifold as a geometric solution of Jang's equation. R. Schoen proposed to employ the observations from \cite{PMTII} to show that closed marginally outer trapped surfaces exist in bounded regions $\Omega \subset M^n$ of the form $\Omega := \Omega_1 - \bar \Omega_2$ for which $\partial \Omega_1$ is outer untrapped and $\partial \Omega_2$ is outer trapped, by using these conditions on the boundary to construct appropriately divergent solutions of Jang's equation on $\Omega$ which would lead to `marginally outer trapped cylinders'
as in \cite{PMTII}. The technical difficulty with this program stems from the one-sidedness of the trapping assumptions, which by themselves are not sufficient to derive a priori boundary gradient estimates to solve the relevant Dirichlet problems for the regularization of Jang's equation. L. Andersson and J. Metzger were able to overcome this difficulty by performing a delicate modification of the data and geometry of $\Omega$ to reduce to the special case where $p \equiv 0$ near $\partial \Omega$, so that the boundary gradient estimate can be derived as for the capillary surface equation. This allowed the authors to prove existence of closed marginally outer trapped surfaces under physically suggestive conditions in \cite{AM}. \\

In this paper we develop a direct and versatile technique to force and control a blow-up of Jang's equation on the original data set, and we use it to prove existence of marginally outer trapped surfaces spanning a given boundary (the Plateau problem). On the way we introduce concepts and low order methods from geometric measure theory to the context of marginally outer trapped surfaces. The robust global geometric features of marginally outer trapped surfaces that are exhibited in this work provide a very tight analogy with the variational aspects of minimal and constant mean curvature surfaces, with immediate applications (cf. Appendix \ref{sec: remark on area of outermost MOTS} and \cite{GAH}). They also open the door to study existence and regularity of marginally outer trapped surfaces in arbitrary dimension.  

\begin{theorem} [Plateau problem for marginally outer trapped surfaces] \label{thm: Plateau problem} Let $(M^n, g, p)$ be
an initial data set and let $\Omega \subset M^n$ be a bounded open domain with smooth boundary $\partial \Omega$. Let $\Gamma^{n-2} \subset \partial \Omega$ be a non-empty, smooth, closed, embedded submanifold that separates this boundary in the sense that $\partial \Omega \setminus \Gamma^{n-2} = \partial_1 \Omega \dot \cup \partial_2 \Omega$ for disjoint, non-empty, and relatively open subsets $\partial_1 \Omega, \partial_2 \Omega$ of $\partial \Omega$. We assume that $\mc_{\partial \Omega} + \tr_{\partial \Omega} p > 0$ near $\partial_1 \Omega$ and that $\mc_{\partial \Omega} -
\tr_{\partial \Omega} p > 0$ near $\partial_2 \Omega$ where the mean curvature scalar is computed as the tangential divergence of the unit normal pointing out of $\Omega$. Then there exists an almost minimizing (relative) boundary $\Sigma^{n-1} \subset \Omega$, homologous to $\partial_1 \Omega$, whose singular set is strictly contained in $\Omega$ and has Hausdorff dimension $\leq (n-8)$,
so that $\Sigma^{n-1}$ satisfies the marginally outer trapped surface equation distributionally, and so that $\Sigma^{n-1}$ is a smooth hypersurface near $\Gamma^{n-2}$ with boundary $\Gamma^{n-2}$. In particular, if $2 \leq n \leq 7$, then $\Sigma^{n-1}$ is a smooth embedded marginally outer trapped surface in $\Omega$ which spans $\Gamma^{n-2}$. \end{theorem}

\noindent The boundary conditions for the domain $\Omega$ in this theorem are physically natural and reduce to a standard assumption (see \cite{HardtSimon}, \cite{MeeksYau}) for the Plateau problem for minimal surfaces in Riemannian manifolds when $p \equiv 0$. We emphasize that $\Sigma^{n-1}$ in this theorem arises as a relative boundary of a subset $E \subset \Omega$ with finite perimeter in $\Omega$ which has the property that \begin{eqnarray*}
|\D \chi_E| (\Omega) \leq |\D \chi_F | (\Omega ) + C \mathcal{L}^n (E \Delta F)   \text{ for every } F \subset \Omega \text{ with } E \Delta F \subset \subset \Omega. 
\end{eqnarray*}
Here $E \Delta F:= (E - F) \cup (F - E)$ is the symmetric difference of $E$ and $F$, $|\D \chi_E|(\Omega)$ and $|\D \chi_F|(\Omega)$ denote the perimeters of $E$ and $F$ in $\Omega$, $\mathcal{L}^n$ is $n$-dimensional Lebesgue measure, and $C := 2 n |p|_{\mathcal{C}(\bar \Omega)}$ is an explicit constant. For convenience we will say throughout that a boundary $\Sigma^{n-1}$ satisfying the above condition is $C$-almost minimizing. This particular property is a special case of much more general notions of almost minimizing currents as pioneered by F. Almgren in \cite{Almgren76}. We review the classical regularity and compactness results for such boundaries, which in essence are those of area minimizing boundaries, and give precise references to the literature in Appendix \ref{sec: almost minimizing property}. Notions of almost minimizing boundaries are typical in the context  of variational problems with obstacles and constraints. We find it surprising that the $C$-almost minimizing property appears quite naturally in our study of marginally outer trapped surfaces, despite the lack of a variational principle. Along with this property, the apparatus of geometric measure theory enters the picture, allowing us to dispense with the use of curvature estimates and to prove existence results in all dimensions. Moreover, area bounds for marginally outer trapped surfaces and the existence and regularity of the apparent horizon which have been established in \cite{AM} in $n=3$ by a delicate surgery procedure based on a priori curvature estimates for stable $2$-dimensional surfaces can be deduced readily from this property and hold in higher dimensions, as we discuss in \cite{GAH}. Besides this `zero order' 
almost minimizing feature, the marginally outer trapped surfaces in the theorem also have certain stability properties (cf. Remark \ref{rem: stronger stability}): if for example the data set satisfies the dominant energy condition, then the first Dirichlet eigenvalue of the conformal Laplacian of $\Sigma^{n-1}$ will be non-negative.  \\

The method presented in this paper immediately adapts to prove
analogous existence results for $\mathcal{C}^{1, \alpha}$-boundaries $\Sigma^{n-1}$ with prescribed boundary manifolds $\Gamma^{n-2}$ that solve $\mc_{\Sigma} = \F(x, \vec
\nu)$ distributionally, provided $\F$ is a continuous function on the Grassmann
manifold $G_{n-1}(\Omega)$. The minor modifications are carried out for the
special case of generalized apparent horizons (i.e. surfaces $\Sigma^{n-1}$
with $\mc_\Sigma = |\tr_{\Sigma}(p)|$, cf. \cite{BrayGAH}) in \cite{GAH}. \\

For context and convenient reference we briefly explain some of the key insights of R. Schoen and S.-T. Yau in their analysis of divergent limits of solutions of the regularized Jang's equation in Proposition 4 of \cite{PMTII}:
\\

\noindent Given a complete asymptotically flat initial data set $(M^n, g, p)$
consider a new ``initial data set'' $(M^n \times \mathbb{R}, g + dx_{n+1}^2,
p)$, where $p$ is extended trivially to $M^n \times \mathbb{R}$ by zero in the
vertical direction. For $u : M^n \to \mathbb{R}$ consider the condition
\begin{equation} \label{eq: Jang's equation}
\left(g^{ij} - \frac{u^i u^j}{1 + |\D u|^2}\right) \left(\frac{\D_i \D_j u}
{\sqrt{1 + |\D u|^2}} + p_{ij} \right) = 0
\end{equation}
This equation was introduced in \cite{Jang} and has become known as
Jang's equation. Note that
\begin{equation*}
\mc(u) := \left(g^{ij} - \frac{u^i u^j}{1 + |\D u|^2}\right) \frac{\D_i \D_j u}
{\sqrt{1 + |\D u|^2}}
\end{equation*}
expresses the mean curvature of the graph of $u$ as a submanifold in $M^n
\times \mathbb{R}$ when computed as the tangential divergence of the downward
pointing unit normal. The expression
\begin{equation*}
\tr(p)(u) := \left(g^{ij} - \frac{u^i u^j}{1 + |\D u|^2}\right) p_{ij}
\end{equation*}
is the trace of $p$ over the tangent space of this graph. Put differently, the graph of a solution $u$ of Jang's equation (\ref{eq: Jang's equation}) is a marginally outer trapped surface in $(M^n \times \mathbb{R}, g + dx_{n+1}^2, p)$. Proving existence of entire solutions $u: M^n \to \mathbb{R}$ of (\ref{eq: Jang's equation}), say with decay to $0$ on the asymptotically flat ends, is hampered by the lack of sub and super solutions for a priori estimates. To overcome this, for parameters $t>0$ the
regularized problems of finding $u_t: M^n \to \mathbb{R}$ with
\begin{equation} \label{eq: regularized Jang's equation}
(\mc + \tr(p)) u_t = t u_t \text{ with }  u_t(x) \to  0 \text{ on the ends of } M^n
\end{equation}
are introduced and studied in \cite{PMTII}. The constant functions $ \pm n |p|_{\mathcal{C}^0} / t$ are barriers for these problems. Together with interior estimates for the gradient, the a priori estimate $t |u_t| + t |\D u_t|  \leq C(|\Ric_M|, |p|_{\mathcal{C}^1})$ is obtained, and existence of solutions $u_t$ of (\ref{eq: regularized Jang's equation}) follows from the continuity method. It is then shown in \cite{PMTII} that vertical translation is an `approximate Jacobi field' on $\graph (u_t)$, where `approximate' is independent of $t>0$. This leads to a uniform stability-like inequality for these graphs. The iteration method from \cite{Schoen-Simon-Yau} provides geometric bounds for these hypersurfaces of $M^n \times \mathbb{R}$ provided that $n \leq 5$. It follows that in low dimensions
one can pass $\graph(u_t)$ to a smooth limiting hypersurface $G \subset M^n \times \mathbb{R}$ along a subsequence $t_i \searrow 0$. This
submanifold $G$ is then a marginally outer trapped surface in $(M^n \times
\mathbb{R}, g+ dx_{n+1}^2, p)$. The Harnack principle implies that every connected component of
$G$ is either itself a graph over some open subset $\Omega_0 \subset
M^n$ (on which the defining function solves Jang's equation), or else
is a cylinder of the form $\Sigma^{n-1} \times \mathbb{R}$ (in which case
$\Sigma^{n-1}$ is a closed marginally outer trapped surface). It also
follows that the graphical components in this limit $G$ must be asymptotic to such marginally outer trapped cylinders
$\Sigma^{n-1} \times \mathbb{R}$ on approach to $\partial \Omega_0$.
By the same token, the boundaries of the sets $\Omega_{+}$ and $\Omega_{-} \subset
M^n$ on which the functions $u_{t_i}(x)$ diverge to plus,
respectively, minus infinity as $t_i \searrow 0$ consist of smooth closed embedded
marginally outer trapped surfaces $\Sigma^{n-1} \subset M^n$. If no such surfaces are
contained in the initial data set or if $|u_t|$ is bounded independently of $t>0$, then this construction leads to an entire solution $u: M^n \to \mathbb{R}$ of Jang's equation.
Importantly, these arguments of \cite{PMTII} show that if solutions $u_t$ of the regularized Jang's equation (\ref{eq: regularized Jang's equation}) diverge to plus infinity in some places and stay
finite or diverge to minus infinity in other places as $t \searrow 0$, then there must be a marginally outer trapped surface $\Sigma^{n-1}$ in the base.  \\

\noindent Note the analogy of this blow-up analysis of \cite{PMTII} with M. Miranda's notion of ``quasi-solutions" of the minimal surface equations (see \S4.1 in \cite{MM}), which provides an appropriate and convenient framework for passing possibly divergent minimal graphs and more generally graphs of functions that solve a variational problem to subsequential limits. \\

R. Schoen suggested to use the features of \cite{PMTII} to prove existence of marginally outer trapped surfaces in domains $\Omega$ by forcing a blow-up of Jang's equation from trapping conditions imposed on the boundary. The existence of closed marginally outer trapped surfaces between an outer trapped and an outer untrapped boundary in \cite{AM} was obtained as a corollary of the analysis in Proposition 4 of \cite{PMTII} explained above by constructing appropriately divergent solutions of certain Dirichlet problems for the regularized Jang's equation on a modified data set. To tackle the Plateau problem from this direction, we are interested in finding solutions $u_t \in
\mathcal{C}^{2, \mu} (\bar{\Omega})$ of \begin{equation} \label{eq: Dirichlet
approximate Plateau}
\begin{array}{rlrrr} (\mc + \tr(p))u_t &= tu_t   & \text{ on } &\Omega& \\
                     u_t        &= \varphi_t & \text{ on }     &\partial \Omega&
\end{array}
\end{equation}
for choices of boundary data $\varphi_t$ that, loosely speaking, diverge to
plus infinity as $t \searrow 0$ on $\partial_1 \Omega$, and to minus infinity on $\partial_2 \Omega$. From \cite{PMTII} one expects that if such solutions $u_t$ exist and if their graphs can be passed to a smooth subsequential limit in $\Omega \times \mathbb{R}$, then this limit either contains a marginally outer trapped cylinder $\Sigma^{n-1}\times \mathbb{R}$, or contains a graphical component asymptotic to such a cylinder, whose cross-section $\Sigma^{n-1} \subset \Omega$ spans the boundary $\Gamma^{n-2} \subset \partial \Omega$. Note that if $(M^n, g)$ is flat and $p \equiv 0$, then (\ref{eq: Dirichlet approximate Plateau}) is just the Dirichlet problem for the ordinary capillary surface equation. In this case the Dirichlet problem is well-posed provided that $t |\varphi_t (x)| \leq \mc_{\partial \Omega} (x)$ at every boundary point $x \in \partial \Omega$, by  \S 20 in \cite{Serrin69}. When $p \not \equiv 0$, the methodology of J. Serrin in \cite{Serrin69} provides sub and super solutions for the boundary gradient estimate for solutions $u_t $ of (\ref{eq: Dirichlet approximate Plateau}) provided that $  - \mc_{\partial \Omega} (x) + \tr_{\partial \Omega} (p)(x) <  t \varphi_t(x) < \mc_{\partial \Omega} (x) + \tr_{\partial \Omega}(p)(x)$ at every boundary point $x \in \Omega$. These barriers are constructed as inward perturbations of the cylinders $\partial \Omega \times \mathbb{R}$ below, respectively above, the boundary values $\varphi_t$. In addition to the one-sided trapping assumption for $\partial \Omega$ at hand this would require the boundary to satisfy the much stronger condition $\mc_{\partial \Omega} > |\tr_{\partial \Omega}(p)|$, or at least to be mean-convex in addition to the original assumption. The original proof of the boundaryless case of the Plateau problem in \cite{AM} is precisely by deforming the  geometry and data carefully to reduce to this special case (in fact to the case where $p \equiv 0$ near $\partial \Omega$). \\

\noindent Our method of constructing appropriately divergent solutions of the regularized Jang's equation takes place on the original data set. However, we will not find exact solutions of the Dirichlet problems. Instead, we are going to use the Perron method, suitably adapted to this context, to obtain maximal interior solutions $u_t$ of the equations $\L_t u_t := (\mc + \tr(p) - t)u_t = 0$ subject to the condition that $\underline{u}_t \leq u_t \leq  \overline{u}_t$. Here, $\underline{u}_t$ and $\overline{u}_t$ are viscosity sub and super solutions for the monotone operators $\L_t$ which coincide with standard barriers obtained from the one-sided trapping assumptions $\mc_{\partial_1 \Omega} + \tr_{\partial_1 \Omega}(p) > 0$ and $\mc_{\partial_2 \Omega} - \tr_{\partial_2 \Omega}(p)  >0$ near the respective portions of $\partial \Omega$, and which are suitably extended to the interior by the constant barriers $\pm n |p|_{\mathcal{C}(\bar \Omega)}/t$. In general these sub and super solutions will only match up in a small neighborhood of $\Gamma^{n-2} \subset \partial \Omega$ (where the trapping assumptions are two-sided), hence forcing $u_t$ to assume prescribed values there. Besides this, we will not be able to say anything precise about the boundary behavior of the $u_t$'s other than to assert their divergence in the regularization limit $t \searrow 0$. Along with the geometric compactness of the graphs of these solutions (which have their mean curvature bounded independently of $t>0$ and hence are almost minimizing) this is fortunately all we need. We note in Lemma \ref{lem: cylindrical almost minimizers} that the almost minimizing property of $\graph(u_t)\subset \Omega \times \mathbb{R}$ descends to the cross-sections of cylinders contained in limits and limits of vertical translations of such surfaces. We can conclude the boundary regularity of $\Sigma^{n-1}$ from an application of \cite{AllardBdry}. \\

\noindent The use of the Perron method in conjunction with the minimal surface operator goes back to the study of Dirichlet problems for minimal and constant mean curvature graphs \cite{JenkinsSerrin66} and \cite{Serrin}: maximal interior solutions to these problems exist provided the domain admits sub and super solutions, and these Perron solutions assume continuous boundary values where the boundary is sufficiently mean-convex so that two-sided barriers can be constructed. The analysis described above shares important features with the classical construction of Scherk-type (i.e. infinite boundary value) minimal graphs in polygonal regions in $\mathbb{R}^2$ in \cite{JenkinsSerrin68} and its obstructions expressed by the Jenkins-Serrin conditions, and variants of it in \cite{Spruck72}, \cite{UM}, and \cite{HRS08}. The elegant reflection method of \cite{Spruck72} to deal with boundary arcs of constant negative curvature does not seem to generalize to our setting. \\  

In \cite{Condensation}, R. Schoen and S.-T. Yau used their work on
Jang's equation to show that if the boundary of a compact domain $\Omega$ is
trapped both ways, i.e. if $\mc_{\partial \Omega} > |\tr_{\partial
\Omega}(p)|$, and if enough matter is concentrated in $\Omega$ compared to its
homological filling radius $\Rad (\Omega)$, then $\Omega$ must contain a marginally outer trapped surface, see also \cite{Yau97}. This important result about the geometry of initial data sets
has been sharpened and generalized by S.-T. Yau in \cite{Yau}. \\

In the next section, we collect several technical lemmas on the regularized Jang's equation. In Section \ref{sec: Perron solutions} we employ the Perron method to construct maximal interior solutions of these equations when only one-sided boundary barriers are available. As an application we present a new proof of the existence of closed marginally outer trapped surfaces between oppositely trapped surfaces, and we record the uniform almost minimizing property of the surfaces that are found. This sets the stage for our main application, the solution of the Plateau problem, in Section \ref{sec: Plateau problem}. Standard results from geometric measure theory regarding the compactness and regularity theory for the classes $\mathcal{F}_C$ of $C$-almost minimizing boundaries, which will be used throughout this paper, are reviewed in Appendix \ref{sec: almost minimizing property}. We conclude with a brief remark on outermost marginally outer trapped surfaces in Appendix
\ref{sec: remark on area of outermost MOTS}.

\begin{acknowledgements*} This work forms part of my thesis and I am very much indebted to my adviser Richard Schoen at Stanford who inspires me in all ways and who suggested this problem to me. I would like to sincerely thank Simon Brendle, Leon Simon, and Brian White for the great example they set and for everything they have taught me over the years. I am grateful to Theodora Bourni for answering my questions regarding some of the statements in \cite{AllardBdry}. Thanks to G. Galloway and Jan Metzger for discussions, questions, and comments.  Finally to Isolde, Sonnenschein: this one is for you.
\end{acknowledgements*}

\section{Interior gradient estimate, the regularized Jang's equation on small balls, and the Harnack principle} \label{sec: interior gradient estimate}

In this section we give three routine lemmas that will be needed later in this paper.  \\

First, we apply the maximum principle method of Korevaar-Simon (see \cite{Korevaar}, \cite{KorevaarSimon}, and \cite{Spruck} for the method on Riemannian manifolds) to obtain interior gradient estimates from one-sided oscillation bounds for solutions of the regularized Jang's equation.  The derivation of the required estimate is well-known and straightforward from these references but not quite contained in their presentation, so we include the proof for completeness. In \cite{Leon76} interior gradient estimates for very general PDE's of prescribed mean curvature type, including those considered in this section, have been obtained by iteration techniques. The author is grateful to Leon Simon for very helpful discussions regarding these estimates. 
\\

\noindent We begin by introducing relevant notation.  Let $(M^n, g, p)$ be an initial data set and let $u$ be a real valued function
defined on an open subset of $M^n$ on which it solves $(\mc + \tr(p))u = t u$ for
some $t \geq 0$. Let $\Sigma^n := \graph (u) \subset M^n \times \mathbb{R}$ and
endow this hypersurface with the induced metric $\bar{g}$ from $(M^n \times \mathbb{R}, g +
dx_{n+1}^2)$. If $x_1, \ldots, x_n$ are local coordinates on $M^n$, we agree on the
following notation and sign convention, cf. \cite{PMTII} and \cite{Spruck}: 
\begin{eqnarray*}
& \bar{g}_{i j} = g_{i j} + u_i u_j & \text{ induced metric in base
coordinates} \\
& \bar{g}^{i j} = g^{i j} - \frac{u^i u^j} {1 + |\D u|^2} & \text{ inverse of
the induced metric}\\
& v = \sqrt{1 + |\D u|^2} & \text{ area element (really area stretch factor)}\\
& \vec{\nu} = \nu^i \frac{\partial}{\partial x_i} -
\frac{1}{v}\frac{\partial}{\partial t}  \text{ where } \nu^i = \frac{u^i}{v} &
\text{ (downward pointing) unit normal } \\
& h_{i j} = v^{-1} \D_{i} \D_j u &\text{ second fundamental form} \\ 
& \mc_{\Sigma} = \mc(u) = \bar{g}^{i j} h_{i j} = \D_i \nu^i &\text{ mean
curvature} \\
& \tr_{\Sigma}(p) = \tr(p)(u) = \bar{g}^{i j} p_{i j} &\text{ trace of } p
\text{ over the tangent space of } \Sigma^n \\
& \Delta_\Sigma u = v^{-1} \mc_{\Sigma}& 
\end{eqnarray*} \noindent Here, all covariant derivatives are taken with respect to the
$g$ metric on the base, and indices are raised with respect to $g$, unless
indicated otherwise. \\

\noindent The following identity for the hypersurface Laplacian of the area element of a graph (see (2.18) in \cite{PMTII} for a derivation) is well-known as the Jacobi
equation (corresponding to the Killing field coming from vertical translation): \begin{equation} \label{eqn: Jacobi equation}
\Delta_{\Sigma} \left( \frac{1}{v} \right)  + \big(\bar{g}(h, h) + \Ric_{M \times \mathbb{R}} (\vec{\nu}, \vec{\nu}) +
\vec{\nu} (\underbrace{t u - \tr(p)(u)}_{\mc_\Sigma}) \big) \frac{1}{v} = 0
\end{equation}
\noindent where $\mc_\Sigma = \mc(u)$ is thought of as a function on $M^n \times
\mathbb{R}$ (independent of the vertical coordinate $x_{n+1}$).  We compute in coordinates on the base that
\begin{eqnarray*}
\vec{\nu} (t u - \tr(p)(u)) &=& \nu^i (t u - \bar{g}^{m n} p_{m n})_{i} \\ &=&
\nu^i (tu - p_{m}^m + \nu^m \nu^n p_{m n})_{i} \\ &=& \frac{t |\D u|^2}{v} -
\nu^i \D_i p_{m}^{m} + 2 p_{m n} h_{i p} \bar{g}^{p n} \nu^i \nu^m  + \D_i p_{m
n} \nu^m \nu^n \nu^i
\end{eqnarray*}
\noindent (all covariant derivatives on $M^n$) so that for a constant $\beta \geq0$
depending only on the size of the Ricci tensor and  $|p|_{\mathcal{C}^1}$ we
have that \begin{eqnarray} \label{estimate for interior gradient estimate}
\bar{g}(h, h) + \Ric_{M \times \mathbb{R}} (\vec{\nu}, \vec{\nu}) + \vec{\nu}
(t u - \tr(p)(u)) \geq \frac{t |\D u|^2}{v} - \beta \geq - \beta.
\end{eqnarray} Note that no assumption on the size of $t|\D u|$ is necessary
here, since it enters the computation with a favorable sign. In view of the ensuing
differential inequality $\Delta_{\Sigma} \left( \frac{1}{v}\right) \leq \frac{\beta}{v}
$, the use of the maximum principle to derive an interior gradient estimate
from a one-sided oscillation bound is now exactly as in \cite{KorevaarSimon}
or \cite{Spruck}:

\begin{lemma} [Interior gradient estimate] \label{lem: interior gradient
estimate}  Fix a point $x_0 \in M^n$, a radius $\rho \in (0, \inj_{x_0} (M^n,
g))$ and $t \geq 0$. Suppose that $u \in \mathcal{C}^3 (\B(x_0, \rho)) \cap
\mathcal{C}(\bar{\B}(x_0, \rho))$ solves $(\mc + \tr(p))u = tu$ on $\B(x_0,
\rho)$ and that for some $T > 0$ either $u(x) \leq u(x_0) + T$ for all $x \in
\B(x_0, \rho)$ or $u(x) \geq u(x_0) - T$ for all $x \in \B(x_0, \rho)$. Then $|\D
u|(x_0)$ is bounded explicitly by an expression in this one-sided oscillation $T$, the
radius $\rho$, and bounds on $\B(x_0, \rho)$ for $|p|$, $|\D p|$, $t|u|$ and the
sectional curvatures of $g$.

\begin{proof} We only consider $u(x) \leq u(x_0) + T$. The other
case follows from this by considering the equation satisfied by $- u$. By decreasing $\rho$ slightly if necessary we can assume that $u \in \mathcal{C}^3(\bar{\B}(x_0, \rho))$. Define
on $\bar{\B}(x_0, \rho)$ the function
$\phi(x):=-u(x_0)+u(x)+\rho-\frac{T+\rho}{\rho^2}\dist(x_0, x)^2$ where $\dist(x_0, x)$
is geodesic distance on $M^n$. Define the open set $\Omega = \{x \in \B(x_0,
\rho) : \phi > 0 \}$ and observe that $\phi = 0$ on the boundary of $\Omega$,
and that $x_0 \in \Omega$. For a number $K \geq 1$ to be chosen later, consider 
the cut-off function $\eta := e^{K \phi} - 1$. At a point
where $\eta v$ assumes its maximum in $\Omega$, one computes that
\begin{eqnarray*}
0 \geq \frac{1}{v} \Delta_\Sigma(\eta v) = \frac{\eta}{v} \Delta_\Sigma v + \frac{2}{v} \bar{g} (\bar \D v, \bar \D \eta) + \Delta_\Sigma \eta \geq \\
\frac{2 \eta}{v^2} |\bar \D v|^2 - \beta \eta + \frac{2} {v} \bar g (\bar \D v, \bar \D \eta) + \Delta_\Sigma \eta = - \beta \eta + \Delta_\Sigma \eta = \\
\beta + \left( K \Delta_{\Sigma} \phi + K^2 |\bar{\D} \phi|^2 -\beta \right) e^{K \phi} \geq \\
K \underbrace{\left( v^{-1} \mc_{\Sigma} - \frac{T + \rho}{\rho^2}
\Delta_{\Sigma} \dist^2(x_0, \cdot) - \beta \right)}_{=: (I)} e^{K \phi} \\
 +        K^2 \underbrace {\left( |\D u - \frac{T + \rho}{\rho^2} \D \dist^2(x_0,
\cdot)|^2 - \frac{(\D u - \frac{T + \rho}{\rho^2} \D \dist^2(x_0, \cdot), \D
u)^2}{1 + |\D u|^2}\right)}_{=: (II)} e^{K \phi}.
\end{eqnarray*}
Recall that $\Delta_{\Sigma} f = \bar{g}^{i j} \D^2_{i j} f - \mc_{\Sigma}
\frac{u^j}{v}\frac{\partial f}{\partial x_j}$ for functions $f$ defined on the
base, so that by the Hessian comparison theorem $\Delta_\Sigma \dist^2(x_0, \cdot)$
is estimated in terms of $\rho$, a bound on the sectional curvatures of $g$ on
$\B(x_0, \rho)$, and a bound on $\mc_\Sigma$ (which in turn is bounded by $n
|p|$ and $t|u|$ from the equation that $u$ satisfies). It follows that $(I)$,
the term multiplying $K$, is a priori bounded, while $\liminf_{|Du| \to \infty}
(II) \geq 1$. By fixing the constant $K$ to be on the order of an absolute
bound for $(I)$, it follows that $|\D u|$ must lie below a definite threshold
at points in $\Omega$ where $\eta v = (e^{K\phi} - 1) \sqrt{1 + |\D u|^2}$
is largest. 
Since the auxiliary function $\phi$ is bounded on
$\Omega$ by $T + \rho$, this implies an upper bound for $\eta v (x_0)$ = $(e^{K
\rho} - 1) \sqrt{1 + |\D u|^2(x_0)}$ and hence for $|\D u|(x_0)$.
\end{proof}
\end{lemma}

\begin{remark} This method of Korevaar-Simon provides interior gradient estimates from one-sided oscillation bounds for $\mathcal{C}^3$-solutions $u$ of $\mc(u) = \F(x, \frac{\D u}{\sqrt{1 + |\D u|^2}}, u)$, provided $\F$ is a $\mathcal{C}^1$-function that is non-decreasing in its last argument.  \end{remark}

In the introduction we discussed how boundary gradient estimates and hence existence for the Dirichlet problem for the regularized Jang's equation (\ref{eq: Dirichlet approximate Plateau}) follow from standard barriers \cite{Serrin} provided the Dirichlet data $\varphi$ and the geometry of $\partial \Omega$ are such that $\mc_{\partial \Omega}(x) > |\tr_{\partial \Omega}(p) - t \varphi|(x)$ at every $x \in \partial \Omega$. In the construction of solutions to the exact Jang's equation (possibly with interior blow-up) with boundary data $\varphi$ this condition needs to be verified only for sufficiently small $t>0$, so the requirement is that $\mc_{\partial \Omega} > |\tr_{\partial \Omega} (p)|$. The sufficiency of this condition was pointed out by R. Schoen and S.-T. Yau. The following general lemma is a straightforward consequence of the continuity method in \cite{PMTII}. The case $t=0$ in this lemma was discussed in \cite{Yau}, \cite{WangYau}, and the case where $\varphi$ is constant and $p \equiv 0$ near  $\partial \Omega$  in \cite{AM}.   

\begin{lemma} \label{lem: Dirichlet for L_t} (Dirichlet problem for regularized Jang's equation) Let $\Omega \subset (M^n, g, p)$ be a bounded domain with smooth boundary, let $t \in (0, 1)$, and assume that $\varphi \in \mathcal{C}(\partial \Omega)$ is such that $\mc_{\partial \Omega} > |\tr_{\partial \Omega} (p) - t \varphi|$.  Then there exists $u \in \mathcal{C}^3(\Omega) \cap \mathcal{C}(\bar \Omega)$ with $(\mc + \tr(p)) (u) = t u $ in $\Omega$ and such that $u = \varphi$ on $\partial \Omega$. 

\begin{proof} First, note that interior extrema of a solution $u$ are a priori bounded by $n |p|_{\mathcal{C}(\bar \Omega)}/t$. As in the proof of Theorem 16.8 in \cite{GT}, successive application of the interior gradient estimates in Lemma \ref{lem: interior gradient estimate}, H\"older estimates for the gradient (\S 13 in \cite{GT}), Schauder theory, and the Arzela-Ascoli theorem readily reduces to the case where $\varphi \in \mathcal{C}^{2, \mu}(\partial \Omega)$. Let $I = \{s \in [0, 1] : \exists w \in \mathcal{C}^{2, \mu} (\bar\Omega) \text{ s.t. } (\mc + s \tr(p)) w = t w \text{ in } \Omega \text{ and } w = s \varphi \text { on } \partial \Omega  \}$. Obviously, $0 \in I$ and the implicit function theorem (see the proof of Lemma 3 in \cite{PMTII} for the linearized operator) shows that $I$ is relatively open. To see that $I$ is closed (and hence all of $[0, 1]$) one requires an a priori estimate for solutions $w \in \mathcal{C}^{2, \mu}$ as in the definition of $I$ that is independent of $s \in I$. By standard elliptic theory (notably the `Bernstein trick' used at the beginning of Section 4 in \cite{PMTII}, H\"older estimates for the gradient as in Theorem 13.8 in \cite{GT}, and Schauder estimates) it is sufficient to derive a boundary gradient estimate for $w$. Observe that by assumption, $\mc_{\partial \Omega} \geq \chi + |s \tr_{\partial \Omega} (p) - t s \varphi|$ for all $s \in [0, 1]$ for some constant $\chi>0$, and that by continuity we may assume that this inequality not only holds for $\partial \Omega$, but also for all surfaces $\partial_\gamma \Omega$ parallel to $\partial \Omega$ at distances $\gamma \in (0, \delta)$ for some $\delta >0$. This is precisely the hypothesis (14.60) of Serrin's boundary gradient estimate as described in Theorem 14.9 in \cite{GT} (with Euclidean metric replaced by $g$): in local coordinates $x_1, \ldots, x_n$ near a point in $\partial \Omega$ the decomposition $a^{ij} = \Lambda a_\infty^{ij} + a^{ij}_0$ and $b = |p| \Lambda b_\infty + b_0$ of (14.43) in \cite{GT} is given by $\Lambda = (1 + |p|^2)$, $a^{ij}_\infty = g^{ij} - \frac{p^i p^j}{|p|^2}$, $a^{ij}_0 = \frac{p^i p^j}{|p|^2}$, $b_\infty = (g^{ij} - \frac{p^i p^j} {|p|^2}) (s p_{ij} - \Gamma_{ij}^k \frac{p_k}{|p|}) - t z$, $b_0 = - t z \frac{\Lambda }{ \Lambda^{1/2} + |p|} - \frac{p^i p^j}{|p|^2} \Gamma_{ij}^k p_k + s p_{ij} \left( g^{ij} \frac{\Lambda}{ \Lambda^{1/2} + |p|} + \frac{p^i p^j}{|p|} \frac{\Lambda^{1/2}}{ \Lambda^{1/2} + |p|} \right)$. The metric $g$ is used to raise indices and compute lengths here. Both $b_0$ and $b_\infty$ are non-increasing in $z$, and the structural conditions (14.59) in \cite{GT} are also clearly satisfied. As for the boundary curvature conditions (14.60), note that $\mathcal{K}^- + b_\infty (x, s \varphi, - \vec \nu) = \mc_{\partial \Omega} + s \tr_{\partial \Omega}(p) - t s \varphi \geq \chi > 0$ (here $\vec \nu$ is the inward unit normal of $\partial \Omega$) leads to a sub solution of the form $\underline{w}= s \varphi - k \dist_{\partial \Omega}$, where $\dist_{\partial \Omega}(x)$ measures geodesic distance to $\partial \Omega$, $\varphi$ is extended to $\Omega$ independently of $\dist_{\partial \Omega}$ through nearest point projection, and where $k \gg 1$ is sufficiently large depending on the (uniform) upper bound for $w$. Analogously, the condition $\mathcal{K}^+ - b_\infty (x, s \varphi, \vec \nu) = \mc_{\partial \Omega} - s \tr_{\partial \Omega} (p) + t s \varphi \geq \chi > 0$ leads to a super solution of the form $\overline{w} = s \varphi + k \dist_{\partial \Omega}$. In fact, if we extend $\vec \nu$ to a neighborhood of $\partial \Omega$ in $\Omega$ as the gradient of the distance function $\dist_{\partial \Omega}$, then $a^{ij}_\infty (x, d \dist_{\partial \Omega}) (\dist_{\partial \Omega})_{ij} + b_{\infty} (x, \overline w, d \dist_{\partial \Omega})  \leq a^{ij}_\infty (x, d \dist_{\partial \Omega}) (\dist_{\partial \Omega})_{ij} + b_{\infty} (x, \varphi (x), d \dist_{\partial \Omega}) = (g^{ij} - \nu^i \nu^j) (\nu_{i;j} + s p_{ij})(x) - t s  \varphi(x) =  \left( - \mc_{\partial_\gamma \Omega} + s \tr_{\partial_\gamma \Omega}(p)  \right)(x) - t s \varphi(x) \leq - \chi $ provided that $\dist_{\partial \Omega}(x) = \gamma < \delta$. This gives the boundary gradient estimate from Theorem 14.9 in \cite{GT}. Note also that by (\ref{eqn: Jacobi equation}), (\ref{estimate for interior gradient estimate}) an interior maximum value of $|\D w|$ in $\Omega$ must be bounded above by $\frac{\beta}{t}$ where $\beta = \beta(s |p|_{\mathcal{C}^1(\bar \Omega)}, |\Ric_M|_{\mathcal{C}(\bar \Omega)})$ (cf. (4.3) in \cite{PMTII}).  \end{proof} 
\end{lemma}

\begin{corollary} \label{cor: r_0} Let $\Omega \subset M^n$ and $t \in (0, 1)$. Fix a constant $C > n |p|_{\mathcal{C}(\bar {\Omega})}$ where $|p|_{\mathcal{C}(\bar \Omega)}$ measures the largest eigenvalue of $p$ with respect to $g$ across $\bar {\Omega}$. For $x \in \Omega$ define $r_D (x) := \frac{1}{2} \min \{ \inj_x (M^n, g), \dist (x, \partial \Omega), r_0(x)  \}$ where $r_0 (x) > 0$ is so small that $\mc_{\partial \B(x, r)} > 3C$ for all $r \in (0, r_0(x))$. Then for every $r \in (0, r_D(x))$ and every $\varphi \in \mathcal{C} (\partial \B(x, r))$ with $|\varphi|_{\mathcal{C}(\partial \B(x, r))} \leq 2C/t$ the Dirichlet problem $(\mc + \tr(p)) u = tu $ in $\B(x, r)$ with $u = \varphi$ on $\partial \B(x, r)$ has a unique solution $u \in \mathcal{C}^{3}(\B(x, r)) \cap \mathcal{C} ( \bar \B(x, r))$. 
\end{corollary}

We will frequently apply the following well-known Harnack argument for limits of graphs. For minimal graphs, this argument was devised in \cite{BombieriGiusti} (see also \S 1 in \cite{LeonRCE}). In the context of Jang's equation a more sophisticated version was derived in Proposition 2 in \cite{PMTII}. 

\begin{lemma} [Harnack principle] \label{Harnack principle} Let $f_k : \Omega \to \mathbb{R}$, $k = 1, 2, \ldots$, be a family of $\mathcal{C}^3$-functions so that for some constant $\beta>0$ it is true that $\Delta_{G_k} \frac{1}{v_k} \leq \frac{\beta}{v_k}$, where $\Delta_{G_k}$ is the (non-positive) Laplace-Beltrami operator of $G_k := \graph(f_k, \Omega) \subset \Omega \times \mathbb{R}$, and where $\frac{1}{v_k} = (1 + |\nabla f_k|^2)^{-\frac{1}{2}}$ is the vertical component of the upward unit normal of $G_k$. If in some open subset $U \subset \Omega \times \mathbb{R}$ such graphs $G_{k}$ converge (in $\mathcal{C}^3$, as normal graphs) to a submanifold $G \subset U$  as $k \nearrow \infty$,  then every connected component of $G$ is either a vertical cylinder (i.e. of the form $\Sigma^{n-1} \times \mathbb{R}$), or is itself a graph over an open subset of $\Omega$. 

\begin{proof}
Observe that $G$ inherits an orientation from $G_k$, and we take its unit normal to be `upward' in this sense. Since $G_k \to G$ in $U$ it follows that $1 / v_k$ converges to the vertical component of the unit normal vector of $G$ (call it $w$). Clearly, $w \geq 0$ and $\Delta_{G} w  - \beta w \leq 0$. By the Hopf maximum principle, either $w$
vanishes identically on each connected component, or else is everywhere positive on it. In the first case, the connected component is evidently cylindrical, in the second alternative it is a graph. 
\end{proof}
\end{lemma}
 
\section{The Perron method for Jang's equation} \label{sec: Perron solutions}

In this section we invoke the Perron method to prove existence of certain maximal interior solutions $u_t$ of the regularized Jang's equations $\L_t u_t := (\mc + \tr(p) - t) u_t = 0$ on domains $\Omega \subset M^n$ whose boundary components satisfy certain one-sided trapping conditions. As a first application, and in preparation for the next section, we supply a new proof of the following theorem proposed by R. Schoen (with existence proven first in \cite{AM} in those dimensions where stability based curvature estimates are available):

\begin{theorem} \label{thm: Schoen's theorem} Let $(M^n, g, p)$ be an initial data set and let $\Omega \subset M^n$ be a connected bounded open subset with smooth embedded boundary $\partial \Omega$. Assume this boundary consists of two non-empty closed hypersurfaces $\partial_1 \Omega$ and $\partial_2 \Omega$ so that
\begin{equation} \label{eq: trapping assumption for Schoen's theorem}
\mc_{\partial_1 \Omega} + \tr_{\partial_1 \Omega} p > 0 \text{ and }
\mc_{\partial_2 \Omega} - \tr_{\partial_2 \Omega} p > 0,
\end{equation} where the mean curvature scalar is computed as the tangential divergence of the unit normal vector field that is pointing out of $\Omega$. Then there exists a closed boundary $\Sigma^{n-1} \subset \Omega$ homologous to $\partial_1 \Omega$ that is $C$-almost minimizing in $\Omega$ for an explicit constant $C = C(|p|_{\mathcal{C}(\bar \Omega)})$ (and hence has singular set of Hausdorff dimension $\leq n-8$) and which satisfies the marginally outer trapped surface equation $\mc_\Sigma + \tr_\Sigma (p) = 0$ distributionally. In particular, if $2 \leq n \leq 7$, then $\Sigma^{n-1}$ is a non-empty smooth embedded closed marginally outer trapped surface in $\Omega$. If $3 \leq n \leq 7$ and if the initial data set satisfies the dominant energy condition, then $\Sigma^{n-1}$ has non-negative Yamabe type.
\end{theorem}

\begin{remark} \label{rem: stronger stability} The fact that the dominant energy condition restricts the Yamabe type of hypersurfaces arising in a blow-up of Jang's equation was recognized and employed in \cite{PMTII}. In \cite{AM}, the authors prove that a closed marginally outer trapped surface arising in the regularization limit of Jang's equation as in \cite{PMTII} satisfies the stronger `non-symmetric' stability condition of \cite{AMS05} that implies the non-negativity of its Yamabe type, as was proven in \cite{GalSch}. \end{remark}

Note that if $\Omega$ arises as $\Omega := \Omega_1 - \bar \Omega_2$ with $\partial_1 \Omega = \partial \Omega_1$ and $\partial_2 \Omega = \partial \Omega_2$, then condition (\ref{eq: trapping assumption for Schoen's theorem}) asks that $\partial \Omega_2$ is outer trapped and $\partial \Omega_1$ is outer untrapped, and the theorem asserts that there exists a closed marginally outer trapped surface $\Sigma^{n-1}$ in between.  \\

As in \cite{AM}, the existence of the marginally outer trapped surface $\Sigma^{n-1}$ in Theorem \ref{thm: Schoen's theorem} here will be an immediate corollary to Proposition 4 in \cite{PMTII}, once appropriately divergent solutions $u_t$ of the regularized Jang's equations $\L_t u_t = 0$ are constructed. R. Schoen pointed out that the trapping assumptions (\ref{eq: trapping assumption for Schoen's theorem}) provide one-sided boundary barriers for Jang's equation but are by themselves not sufficient to solve the relevant Dirichlet problems (\ref{eq: Dirichlet approximate Plateau}). However, these one-sided boundary barriers together with the constant sub and super solutions coming from the capillarity regularization of Jang's equation do lead to \emph{maximal interior solutions} $u_t$ of $\L_t u_t = 0$ which have the \emph{required blow-up} behavior. The goal of this section is to make this statement precise and to prove it. Unlike with more typical Dirichlet problems, the solutions $u_t$ here will not in general assume particular boundary values anywhere. Importantly though, we retain uniform and explicit control on $|t u_t|$ in our construction and hence on the mean curvature of $\graph(u_t) \subset M^n \times \mathbb{R}$, so that we are able to pass these hypersurfaces to a subsequential limit in the class of $C$-almost minimizing boundaries. This also means that we can dispense with the curvature estimates used in \cite{PMTII} and \cite{AM} which are limited to low dimensions. \\

\noindent The original proof of Theorem \ref{thm: Schoen's theorem} in \cite{AM} approaches the construction of `solutions' 
$u_t$ of the approximate Jang's equation by reduction to a special case. The authors first bend $\Omega$ along its boundary towards its `lightcone' 
in the Lorentzian manifold $(M^n \times (-\varepsilon, \varepsilon), - dt^2 + g + tp)$ (where $\varepsilon > 0$ is small), carefully, so that the altered data has everywhere positive mean curvature $\tilde \mc_{\partial \Omega} > \delta$, and so that (after a further modification through a cut-off) $\tilde p$ vanishes near the boundary. After this change of the boundary geometry and data, the regularized Jang's equations coincide with capillary surface equations in a neighborhood of $\partial \Omega$. At this point the authors can use standard, two-sided boundary barriers for boundary values $\pm \frac{\delta}{2t}$, arising as in Theorem 14.9 in \cite{GT} from linear inward perturbations of the cylinder $\partial \Omega \times \mathbb{R}$ above and below the boundary values, to obtain the a priori estimates for solving the relevant Dirichlet problems classically. (The same linear barriers apply under one-sided trapping assumptions in the general case, and we will use them in Lemma \ref{lem: toy sub and super solutions}.) The authors then apply the maximum principle to argue that the marginally outer trapped surface they construct from this modified data actually lies in the unchanged part of $\Omega$, and hence is a solution to the original problem. \\

Our construction of appropriately divergent solutions $u_t$ of the regularized Jang's equation will occupy the rest of this section and will at the same time set the stage for the
more delicate existence results of Section \ref{sec: Plateau problem}. To set up the appropriate definition of Perron sub and super solutions, we fix a constant $C > n |p|_{\mathcal{C}(\bar \Omega)}$ and let $r_D: \Omega \to \mathbb{R}^+$ be the function from Corollary \ref{cor: r_0}. The subscript $D$ here will remind us that the Dirichlet
problem $\L_t u_t = 0$ on $\B(x, r)$ for the regularized Jang's equation can be solved for arbitrary continuous boundary data $\varphi$ on $\partial \B(x, r)$ with $|\varphi| \leq \frac{2C}{t}$, for any $x \in \Omega$ and $r \in (0, r_D(x))$. For convenience we will adjust the definition of $r_D$ slightly just before the statement of Lemma \ref{lem: toy sub and super solutions}. With this preparation, the definition of Perron sub and super solutions and their elementary properties in the next few lemmas are standard. We refer to \S 2.8 in \cite{GT} for the classical case of harmonic functions, and to \cite{JenkinsSerrin66}, \cite{JenkinsSerrin68}, \cite{Serrin}, \cite{Spruck72}, \cite{HRS08} for variants of the Perron method in conjunction with Dirichlet problems for constant mean curvature graphs.  \\

\begin{definition} \label{defn: Perron sub and super solution} A Perron sub
solution of $\L_t$ on $\Omega$ is a function $\underline{u} \in
\mathcal{C}(\bar{\Omega})$  with $|\underline{u}| \leq \frac{C}{t}$ and the
following property: if $x \in \Omega$, $r \in (0, r_D(x))$ and $u \in
\mathcal{C}(\bar{\B}(x, r)) \cap \mathcal{C}^{2, \mu}(\B(x, r))$ is the unique
solution of the Dirichlet problem \begin{eqnarray*} 
\begin{array}{rlllr} \L_t u &= & 0              &\text{ on } &\B(x, r) \\
                        u  &= & \underline{u}  &\text{ on } &\partial \B(x, r)
\end{array}
\end{eqnarray*}
whose existence is guaranteed by Theorem \ref{lem: Dirichlet for L_t}, then $u
\geq \underline{u}$ on $\bar{\B}(x, r)$. A function $\bar{u} \in \mathcal{C}(\bar \Omega)$  with $|\bar u| \leq \frac{C}{t}$ is a Perron super solution if for every $x \in \Omega$ and every $r \in (0, r_D(x))$ the unique solution $u$ of the Dirichlet problem 
\begin{eqnarray*}
\begin{array}{rlllr} \L_t u &= & 0              &\text{ on } &\B(x, r) \\
                        u  &= & \bar{u}  &\text{ on } &\partial \B(x, r)
\end{array}
\end{eqnarray*}
lies below $\bar u$ on $\B(x, r)$. Observe that by the maximum principle and the choice of $C \geq n |p|_{\mathcal{C} (\bar \Omega)}$, the constant functions $ - \frac{C}{t}$ and $\frac{C}{t}$ are, respectively, Perron sub and super solutions  of $\L_t$ on $\Omega$. 
\end{definition}

\begin{definition} [$\L_t$ lift of a Perron sub solution, cf. \S 2.8 in \cite{GT}] Let $\underline{u}$
be a Perron sub solution of $\L_t$ on $\Omega$, and fix $x\in\Omega$ and $r \in
(0, r_D(x))$. Let $u \in \mathcal{C}(\bar{\B}(x, r)) \cap \mathcal{C}^{2,
\mu}(\B(x, r))$ be the unique solution of $\L_t u = 0 $ on $\B(x, r)$ with $u =
\underline u$ on $\partial \B(x, r)$. The function $\underline {\hat u}$ defined
by

\begin{equation*} \underline {\hat u} : = \left \{ \begin{array}{lrr} u &
\text{ on } & \B(x, r) \\
                           \underline u & \text{ on } & \bar \Omega \setminus
\B(x, r) \end{array} \right.
\end{equation*}
is called the $\L_t$ lift of $\underline u$ with respect to $\B(x, r)$.

\end{definition}

\begin{definition} \label{defn: class S} Let $\bar{u}_t$ be a Perron super
solution of $\L_t$ on $\Omega$. We denote by $\mathcal{S}_{\bar{u}_t} := \{
\underline{u} \in \mathcal{C}(\bar\Omega) : \underline{u} \text{ is a sub
solution of } \L_t  \text{ on } \Omega \text{ with } \underline{u} \leq
\bar{u}_t \}$ the class of all Perron sub solutions of $\L_t$ lying below
$\bar{u}_t$.
\end{definition}

\begin{lemma} [Basic properties of Perron sub solutions, cf. \S2.8 in \cite{GT}] \label{lem: properties
of sub solutions} Let $\bar u_t$ be a Perron super solution of $\L_t$ on
$\Omega$. We have the following two basic properties of the class $\mathcal{S}_{\bar{u}_t}$:
\begin{enumerate}
\item If $\underline u, \underline v \in \mathcal{S}_{\bar{u}_t}$, then
$\max\{\underline u, \underline v\} \in \mathcal{S}_{\bar{u}_t}$.
\item If $\underline u \in \mathcal{S}_{\bar u_t}, x \in \Omega, r \in (0,
r_D(x))$ and $\underline {\hat u}$ is the $\L_t$ lift of $\underline u$ with
respect to $\B(x, r)$, then $\underline {\hat u} \in \mathcal{S}_{\bar{u}_t}$.
\end{enumerate}
\begin{proof} If $t>0$, $U$ is a smooth bounded domain, and if $f, g \in \mathcal{C}^2(U) \cap \mathcal{C}(\bar {U})$ satisfy $\L_t f \leq \L_t g$ and $f \geq g $ on $\partial U$, then $f \geq g$ on $\bar {U}$, by the maximum principle. Both properties are immediate consequences of this.  \end{proof}
\end{lemma}

\begin{definition} [cf. \S 2.8 in \cite{GT}] \label{defn: Perron solution} Let $\bar u_t$ be a Perron
super solution of $\L_t$ on $\Omega$ such that $\mathcal{S}_{\bar u _t} \neq
\emptyset$. Define the Perron solution $u^P_t$ of $\L_t$ on $\Omega$ with
respect to $\bar u_t$ pointwise for $x \in \Omega$ by $u^P_t (x) := \sup \{
\underline u (x) : \underline u \in \mathcal{S}_{\bar u _t} \}$.
\end{definition}

The interior gradient estimates from Section \ref{sec: interior gradient estimate} and the resulting compactness of bounded solutions (cf. \cite{Serrin} and Corollary 16.7 in \cite{GT}) now imply the interior regularity of the Perron solution in a standard way:

\begin{lemma} [cf. \S 2.8 in \cite{GT} and \cite{Serrin}] \label{lem: existence of Perron solution} Let $u^P_t$ be the Perron solution constructed under the assumptions of Definition \ref{defn: Perron solution}. Then $u^P_t \in \mathcal{C}^{2, \mu}_{\text{loc}}(\Omega)$ and $\L_t u^P_t = 0$ holds classically.
\begin{proof} We will make repeated use of the properties of $\mathcal{S}_{\bar u _t}$ stated in Lemma \ref{lem: properties of sub solutions}. Fix $x \in \Omega$, $r \in (0, r_D(x))$, and let $\underline u_i \in \mathcal{S}_{\bar u _t}$ be a sequence of sub solutions such that $\lim_{i \to \infty} \underline u_i(x) = u^P_t(x)$. By Lemma \ref{lem: properties of sub solutions} we may assume that $\underline u_i$ are pointwise non-decreasing on $\Omega$ and that they solve $\L_t \underline u_i = 0$ classically on $\B(x, r)$. Since $|\underline u_i| \leq \frac{C}{t}$ on $\Omega$ we can use the interior gradient estimate in  Lemma \ref{lem: interior gradient estimate} together with H\"older estimates for the gradient (\S 13 in \cite{GT}) and Schauder estimates to show that on $\B(x, r)$ the limit $\lim_{i \to \infty} \underline u_i(y) =: \tilde u^P_t (y)$ exists in $\mathcal{C}^{2, \mu}_{\text{loc}}$ and that $\L_t \tilde u^P_t = 0$ there. By our choices and the definition of the Perron solution $u^P_t$ we also see that $\tilde u^P_t(y) \leq u^P_t(y)$ for all $y \in \B(x, r)$ with equality at $x$. If for some $y \in \B(x, r)$ we had $\tilde u^P_t (y) < u^P_t(y)$ there would exist a pointwise non-decreasing sequence of functions $\underline v_i \in \mathcal{S}_{\bar u _t}$ with $\underline v_i \geq \underline u_i$ on $\Omega$ such that $\L_t \underline v_i = 0$ on $\B(x, r)$ and with $\lim_{i \to \infty} \underline v_i(y) = u^P_t(y) > \tilde u^P_t(y)$. Applying the same reasoning as above we could then pass $\underline v_i$ on $\B(x, r)$ to a $\mathcal{C}^{2, \mu}_{\text{loc}}$ limit $v$ with $\L_t v = 0$ and $v \geq \tilde u_t^P$ with equality at $x$ and strict inequality at $y$. Since $v - \tilde u^P_t \geq 0$ solves a linear elliptic equation on $\B(x, r)$ this is in contradiction with the strong maximum principle. Hence $\tilde u^P_t = u^P_t$ on $\B(x, r)$ so that $u^P_t \in \mathcal{C}^{2, \mu}_{\text{loc}}$ is a classical solution on $\Omega$, as asserted. \end{proof}
\end{lemma}

Recall that we chose $C > n |p|_{\mathcal{C}(\bar \Omega)}$. We now fix constants $\chi \in (0, C)$ and  $\delta \in (0, \frac{1}{5} \dist(\partial_1 \Omega, \partial_2 \Omega))$ so that the distance surfaces $\partial_{1, \gamma} \Omega := \{ x \in \Omega : \dist(x, \partial_1 \Omega) = \gamma\}$ and analogously $\partial_{2, \gamma} \Omega$ are smooth, disjoint, and embedded for all $\gamma \in [0, 2 \delta]$, and so that  $\mc_{\partial_{1, \gamma} \Omega} + \tr_{\partial_{1, \gamma} \Omega}(p) \geq 2 \chi$ and $\mc_{\partial_{2, \gamma} \Omega} - \tr_{\partial_{2, \gamma} \Omega}(p) \geq 2 \chi$. The methodology of J. Serrin \cite{Serrin69}  applies here and relates the trapping assumptions on the boundary to the existence of sub and super solutions for the operator $\L_t$ in a natural and geometric way. In particular, the proof of Theorem 14.9 in \cite{GT} (or the computation in \S5 of \cite{Yau}) shows that $(\mc + \tr(p))(- k \dist_{\partial_1 \Omega}(x)) \geq \chi$ on $\Omega' := \{x \in \bar \Omega : \dist(x, \partial_1 \Omega) \leq 2 \delta \}$ provided $k$ is sufficiently large. By translation invariance it follows that the function $\underline{u}'_t (x):= \frac{\chi}{t} - k \dist_{\partial_1 \Omega} (x)$ with $k := \frac{\chi + C}{\delta t}$ satisfies $\L_t \underline{u}_t' \geq 0$ on $ \Omega' $ provided $t >0$ is sufficiently small. As in \cite{AM} we choose the boundary value of $\underline{u}_t'$ on $\partial_1 \Omega$ essentially as large as condition (14.60) in \cite{GT} warrants, so as to achieve a blow-up as $t \searrow 0$. Geometrically, the graph of $\underline{u}_t'$ is a linear inward perturbation of the boundary cylinder $\partial_1 \Omega \times \mathbb{R}$ below height $\frac{\chi}{t}$.  We briefly pursue this geometric point of view to elucidate and justify the choice of $\underline{u}_t'$ here: \\

\noindent First, our choice of the outward unit normal for the computation of the mean curvature of $\partial_1 \Omega \times \mathbb{R}$ corresponds to the downward unit normal of the graph of $\underline{u}_t'$. If $\gamma \in [0, 2 \delta]$ and if $N^n = \partial_{1, \gamma} \Omega \times \mathbb{R}$, then $\mc_N + \tr_N (p) - t x_{n+1}$, as a function on $N^n$, is $\geq \chi$ on the part of $N^n$ that lies below $\{ x_{n+1} = \frac{\chi}{t} \}$. More generally, consider the local geometric operator  $\mathcal{L}_t : N^n \to \mc_N + \tr_N(p) - t x_{n+1}$ which takes an oriented hypersurface $N^n$ to a scalar functions on $N^n$. If $\tilde{N}^n$ is a normal graph above some hypersurface $N^n$ with defining function $\phi$, and if we identify functions on $N^n$ and $\tilde{N}^n$ via this parametrization, then $|\mathcal{L}_t N - \mathcal{L}_t \tilde N| \leq c(\{\spt (\phi) \subset N\}) |\phi|_{\mathcal{C}^2(N)}$. Now $\mathcal{L}_t (\graph (\underline{u}'_t, \Omega') ) (x, \underline{u}'_t(x)) = (\L_t \underline {u}_t') (x)$ at every point $x \in \Omega'$. Note that if $\dist(x, \partial_1 \Omega) = \gamma$ and if $t>0$ is sufficiently small, then $\graph( \underline{u}'_t, \Omega')$ can be expressed as a normal graph over the cylinder $\partial_{1, \gamma} \Omega \times \mathbb{R}$ on a ball of radius $1$ with center at $(x, \underline{u}_t'(x))$ in $M^n \times \mathbb{R}$, and that the $\mathcal{C}^2$-norm of the defining function is on the order of $t$. Since $\graph( \underline{u}'_t, \Omega')$ lies below $\{ x_{n+1} = \frac{\chi}{t} \}$, it follows that if $t>0$ is sufficiently small then certainly $ \L_t \underline {u}_t' (x) \geq 0$, as needed. \\

\noindent Finally, we let $\underline{u}_t := \min\{ \underline{u}'_t, - \frac{C}{t} \}$ on $\Omega'$. Note that this extends to a continuous function on all of $\bar \Omega$ with constant value $- \frac{C}{t}$ away from the boundary. At this point, it is convenient to re-define $r_D(x)$ from Corollary \ref{cor: r_0} slightly as $\min\{ r_D(x), \delta\}$. To see that $\underline{u}_t$ is a Perron sub solution, note that the condition in Definition \ref{defn: Perron sub and super solution} only needs to be verified on small balls $\B(x, r)$ where $r \in (0, r_D(x))$. It follows that on any such ball $\B(x, r)$, either $\underline{u}_t \equiv -\frac{C}{t}$ or else $\B(x, r) \subset \Omega'$. From this it follows as in Lemma \ref{lem: properties of sub solutions} that $\underline{u}_t$ is a Perron sub solution. The argument for the construction of the super solution $\bar{u}_t$ is completely analogous. In conclusion we have the following lemma, which provides divergent one-sided Perron sub and super solutions from one-sided trapping assumptions on the boundary (cf. the proofs of Theorem 14.9 in \cite{GT} and of Proposition 3.5 and Lemma 3.7 in \cite{AM}). \\

\begin{lemma} [Existence of special Perron sub and super solutions] \label{lem:
toy sub and super solutions} For all $t>0$ sufficiently small the following
$\mathcal{C}(\bar\Omega)$ functions
\begin{equation*}
\underline u_t (x) := \left\{ \begin{array} {rl} \frac{\chi}{t} -
\frac{\chi + C}{t \delta} \dist_{\partial \Omega}(x) & \text  { if }
\dist (\partial_1 \Omega, x) \leq \delta  \\ - \frac{C}{t} & \text{ if }
\dist (\partial_1 \Omega,  x) > \delta \end{array}  \right.
\end{equation*}
\noindent and
\begin{equation*}
\bar u_t (x) := \left\{ \begin{array} {rl} - \frac{\chi}{t} +
\frac{\chi + C}{t \delta} \dist_{\partial \Omega}(x) & \text  { if }
\dist (\partial_2 \Omega, x) \leq \delta  \\ + \frac{C}{t} & \text{ if }
\dist (\partial_2 \Omega, x) > \delta \end{array}  \right.
\end{equation*}
\noindent are, respectively, Perron sub and super solutions of $\L_t$ on
$\Omega$. We have that $\underline u_t < \bar u_t$ on $\Omega$ and that
$\underline u_t \geq \frac {\chi}{2t}$ on the collar of width
$\frac{\chi \delta}{2(\chi + C)}$ about $\partial_1\Omega$ and that
$\bar{u}_t \leq - \frac{\chi}{2t}$ on the collar of width $\frac{\chi
\delta}{2(\chi + C)}$ about $\partial_2 \Omega$.
\end{lemma}

We now finish the proof of Theorem \ref{thm: Schoen's theorem}. For all
$t>0$ small enough for Lemma \ref{lem: toy sub and super solutions} to apply,
we let $u_t^P \in \mathcal{C}^{2, \mu}_{\text{loc}} (\Omega)$ be the Perron solution $\L_t u_t^P = 0$ from Lemma \ref{lem: existence of Perron solution} for the (non-empty) class $\mathcal{S}_{\bar
u_t}$. The relative boundary of the region $\{(x, x_{n+1}) : x_{n+1} > u_t(x)\}$ in $\Omega \times \mathbb{R}$ is the graph $G_t := \graph(u_t^P, \Omega) \subset \Omega \times \mathbb{R}$ whose mean curvature is clearly bounded by $2C$. Restricting
to the open region $U = \{x \in \Omega: \dist(x, \partial \Omega) > \frac{\chi \delta}{3(\chi + C)} \} \times \mathbb{R}$, we can use the compactness and regularity theory discussed in Appendix \ref{sec: almost minimizing property} to show that in dimensions $2 \leq n \leq 6$, the graphs $G_t$ pass to a smooth subsequential limit $G$ in the class $\mathcal{F}_{2C}$ as $t \searrow 0$, and that this limit $G$ has distance at least $\frac{\chi \delta}{2(\chi + C)}$ from $\partial \Omega \times \mathbb R$. We explain in Remark \ref{rem: better regularity for graphs} why we get smooth convergence also if $n = 7$. If $n \geq 8$, this convergence is in the sense of Lemma \ref{lem: in F_C, current convergence is varifold convergence} (i.e. current \emph{and} varifold convergence) and $\dim \sing G \leq n-8$. The limit $G$ is $2C$-almost minimizing. Note that $G$ is non-empty: the mean-value theorem and the monotonicity formula show that $G$ must have non-empty intersection with every horizontal plane. The remainder of the blow-up analysis is now identical to Proposition 4 in \cite{PMTII}: if $G_0$ is a connected component of the regular set of $G$, then there is a canonical choice of orientation on $G_0$ (the `downwards' 
orientation that $G_0$ inherits from being a geometric limit of graphs; note that Allard's theorem applies to ensure smooth convergence towards $G_0$) so that $\mc_{G_0} + \tr_{G_0} (p) = 0$ holds. This holds distributionally on all of $G$ because the singular set of $G$ is small. By the Harnack principle for limits of graphs in Lemma \ref{Harnack principle}, each $G_0$ is either cylindrical or a graph of a smooth function $u : \Omega_0 \to \mathbb{R}$ defined on some open subset $\Omega_0 \subset \Omega$ where $u$ solves Jang's equation $(\mc + \tr(p))u = 0$. The same discussion applies to any (necessarily non-vanishing) subsequential limit in $\mathcal{F}_{2C}$ obtained from $G$ by a divergent sequence of downward translations.  Such a limit is bound to be a marginally outer trapped cylinder $\Sigma^{n-1} \times \mathbb{R}$ where $\Sigma^{n-1}$ is a relative boundary in $\Omega$ homologous to $\partial_1 \Omega$. By Lemma \ref{lem: cylindrical almost minimizers} the cross-section $\Sigma^{n-1}$ is a $2C$-almost minimizing boundary, and clearly it is itself marginally outer trapped. \\

\noindent Finally, the statement about the Yamabe type of $\Sigma^{n-1} \subset \Omega$ in dimensions $3 \leq n \leq 7$ follows, exactly as in \cite{PMTII}, from the inequality 
\begin{eqnarray*}
0 \leq \bar{\R} - \sum_{i , j = 1}^n (h_{ij} + p_{ij})^2 - 2 \sum_{i=1}^n (h_{i n+1} +
p_{i n+1})^2 + 2 \sum_{i=1}^n \bar{\D}_i (h_{i n+1} + p_{i n+1}) \\ \nonumber +
(\sum_{i=1}^n p_{ii})^2 - \mc^2 +  2 \big( \sum_{i=1}^n p_{ii} + \mc
\big) p_{n+1 n+1} - 2 e_{n+1} \big( \sum_{i=1}^n p_{ii} + \mc \big)
\end{eqnarray*}
(derived as (2.15) in \cite{PMTII}) which is a consequence of the Jacobi equation (\ref{eqn: Jacobi equation}), the Codazzi and Gauss equations, and the dominant energy condition. The indices $\{i, j\}$ here are with respect to a local orthonormal frame $\{e_1, \ldots , e_{n+1}\}$ for $\Omega \times \mathbb{R}$ which is parallel in the vertical direction and such that $e_{n+1}$ restricts to the downward pointing unit normal field on $G_t$. The function $\mc$ on $G_t$ is also extended vertically here; $\bar{\R}$ and $\bar{\D}$ are used to denote the scalar curvature and covariant differentiation on $G_t$ (making the third term on the right the divergence of the tangential vector field whose length is computed in the second term), $h_{ij}$ is the second fundamental form of $G_t$, and $h_{i n+1} = - \bar{\D}_i \log{\sqrt{1 + |\D u_t|^2}}$. We point out that it is important here that $e_{n+1} \big( \sum_i p_{ii} + \mc \big) = e_{n+1} (tu_t) = \frac{t |\D u_t|^2}{\sqrt{1 + |\D u_t|^2}}$ enters the inequality with a favorable sign. 

\begin{remark} \label{rem: blow up is n |p| almost minimizing} In the smooth case $2 \leq n \leq 7$, the marginally trapped cylinder $\Sigma^{n-1} \times \mathbb{R}$, obtained above as a limit of graphs of solutions $u_{t_i}$ of the regularized Jang's equations with $t_i |u_{t_i}| \leq C$ followed by a sequence of divergent downward translations, is always $n |p|_{\mathcal{C}(\bar \Omega)}$-almost minimizing. To see this, consider the relatively open sets $\Omega_+, \Omega_-, \Omega_0$ where, respectively, the $u_{t_i}$ either diverge to $+ \infty$, to $- \infty$, or converge. By the analysis in Proposition 4 in \cite{PMTII}, the whole domain $\bar \Omega$ is the union of the closures of these sets. Recall that $ \mc(u_{t_i}) $ is the divergence of the unit normal vector field of $\graph(u_{t_i})$ that leads to the almost minimizing property in Example \ref{ex: graphs of bounded mean curvature}. If $\Omega'_+ \subset \subset \Omega_+$ then $\liminf_{i \to \infty} \inf_{x \in \Omega'_+} \mc(u_{t_i}) (x) \geq - n |p|_{\mathcal{C}(\bar \Omega)}$, and similar statements hold for the other regions. A careful analysis of signs in Example \ref{ex: graphs of bounded mean curvature} then establishes the asserted $n |p|_{\mathcal{C}(\bar \Omega)}$-almost minimizing property of $\Sigma^{n-1} \subset (\partial \Omega_+ \cup \partial \Omega_- \cup \partial \Omega_0) \cap \Omega$.
\end{remark}

\section{The Plateau problem for marginally outer trapped surfaces} \label{sec: Plateau
problem}

In this section we establish the existence of marginally outer trapped surfaces spanning a given boundary (Theorem \ref{thm: Plateau problem}). For the duration of the proof we fix a constant $C > n |p|_{\mathcal{C}(\bar \Omega)}$. Let $\dist_\Gamma$ measure oriented geodesic distance in $\partial \Omega$ to $\Gamma^{n-2}$ (say positive
on $\partial_1 \Omega$). Let $\delta > 0$ and $\chi \in (0, C)$ be positive
constants such that $\dist_\Gamma$ is smooth on $\{x \in \partial \Omega : - 2 \delta < \dist_\Gamma (x) < 2 \delta \}$ and such that \begin{eqnarray*} 
& \mc_{\partial \Omega} + \tr_{\partial \Omega} p \geq 2 \chi & \text{ on  }
\partial_1 \Omega \cup \{x \in \partial \Omega: - 2 \delta < \dist_\Gamma(x) < 2
\delta \} ,  \\ \nonumber &\mc_{\partial \Omega} - \tr_{\partial
\Omega} p \geq 2 \chi & \text{ on } \partial_2 \Omega \cup \{x \in \partial
\Omega : - 2 \delta < \dist_\Gamma(x) < 2 \delta \}.
\end{eqnarray*} 

\begin{lemma} \label{lem: barriers for Plateau problem} For every $t \in (0, 1)$ there
exists a Perron sub solution $\underline{u}_t \in \mathcal{C}(\bar \Omega)$ of
$\L_t$ on $\Omega$ (in the sense of Definition \ref{defn: Perron sub and super
solution}) such that for $x \in \partial \Omega$ we have
\begin{displaymath}
\underline{u}_t (x) = \text{smooth version of} \left\{ \begin{array}
{rlrll}
- \frac{C}{t}                                  &             &&\dist_\Gamma(x)& <
- \delta t \\
  \frac{C}{\delta t^2} \dist_\Gamma (x)              & \text{ if } &- \delta t <
&\dist_\Gamma(x)& < \frac{\chi \delta t}{C}    \\
  \frac{\chi}{t}                           &             &\frac{\chi
\delta t }{C} < &\dist_\Gamma(x).
\end{array} \right.
\end{displaymath}

\noindent Analogously there exists a Perron super solution $\bar{u}_t \in
\mathcal{C}(\bar \Omega)$ with

\begin{displaymath}
\bar {u}_t (x) = \text{smooth version of} \left\{ \begin{array} {rlrll}
- \frac{\chi}{t}                                  &
&&\dist_\Gamma(x)& < - \frac{\chi \delta t}{C} \\
  \frac{C}{\delta t^2} \dist_\Gamma (x)              & \text{ if }         &-
\frac{\chi \delta t}{C} < &\dist_\Gamma(x)& < \delta  t   \\
  \frac{C}{t}                           &             &\delta t <
&\dist_\Gamma(x).
\end{array} \right.
\end{displaymath}

\noindent We can arrange so that $- \frac{C}{t} \leq \underline{u}_t \leq
\bar{u}_t \leq \frac{C}{t}$ on $\bar{\Omega}$.

\begin{proof}
The sub solution $\underline{u}_t$ can be obtained in a standard way \cite{GT} from linear inward perturbation of $\partial_1 \Omega \times \mathbb{R}$ beneath the required boundary values. We argue this geometrically here. First, note that from the trapping assumptions we know that $\mathcal{L}_t (\partial \Omega \times \mathbb{R}) \geq \chi$ holds on $(U \cap \bar \Omega) \times (- \infty, \frac{\chi}{t})$, where $U$ is an open neighborhood of $\partial_1 \Omega \cup \{x \in \partial \Omega : -  \delta < \dist_\Gamma (x) <  \delta \}$ in $M^n$. $\mathcal{L}_t$ is the geometric regularized Jang's operator that was introduced in the proof of Lemma \ref{lem: toy sub and super solutions}. Let $h_{\partial \Omega}: \partial \Omega \to \mathbb{R}$ denote the intended boundary values of the sought-after sub solution $\underline{u}_t$. Extend the inward pointing unit normal $\vec \nu (x)$ of $\partial \Omega$ to a normal vector field of $\partial \Omega \times \mathbb{R}$ by vertical translation. For large constants $k \gg 1$ consider the function $d : \partial \Omega \times \mathbb{R} \to \mathbb{R}$ defined by $d (x, x_{n+1}) = \frac{1}{k} (h_{\partial \Omega}(x) - x_{n+1})$. For $k$ sufficiently large (depending on $t>0$) the normal graph $N^n : = \{ \exp_{(x, x_{n+1})}(d(x, x_{n+1}) \vec \nu(x)) :    (x, x_{n+1}) \in \partial \Omega \times (-\frac{C+1}{t}, \frac{\chi}{t}) \} \cap ((U\cap \bar \Omega) \times \mathbb{R})$ coincides with the vertical graph of a function $\underline{u}_t'$ that is defined on a subset $U' \subset U\cap\bar\Omega$, because $d$ is monotone in $x_{n+1}$. Since $|d|_{\mathcal{C}^2} = O (\frac{1}{k})$, it follows that $\mathcal{L}_t (N^n) = \mc_N + \tr_N (p) - t x_{n+1}  \geq 0 $ on $U' \times \mathbb{R}$ if $k$ is sufficiently large (again depending on $t>0$). Hence $\L_t (\underline{u}_t')  \geq 0$ on $U'$. By increasing $k$ further if necessary we can arrange for $\underline{u}_t := \max (- \frac{C}{t}, \underline u_t')$ to extend to a continuous function on all of $\bar \Omega$. This will be the required sub solution. The construction of the super solution is identical.
\end{proof}
\end{lemma}

\noindent The sub and super solutions of the preceding lemma are not optimal
for our purposes away from $\Gamma^{n-2}$ in that they wouldn't prevent the
graphs of the corresponding Perron solutions $u_t^P$ from converging towards
$(\partial \Omega \setminus \Gamma^{n-2}) \times \mathbb{R}$ at fixed bounded
heights as $t \searrow 0$. In order to get this shielding effect from the
initial barriers we could improve the construction in Lemma \ref{lem: barriers
for Plateau problem} by allowing the inward slope (respectively the normal
slope $k$ of the perturbation) of $\bar{u}_t$ and $\underline{u}_t$ to vary as a
function of $\dist_\Gamma(\cdot)$. Alternatively, one can construct for every $x \in
\partial_1 \Omega$ a local sub solution whose graph is a cone with vertex at
$(x, \frac{\chi}{t})$ and base a small inward perturbation (depending on
$\dist_\Gamma(x)$ but not on $t>0$) of $\partial \Omega \times \{- \frac{C}{t} \}$ near $(x, -\frac{C}{t})$, continued by
the horizontal plane at height $- \frac{C}{t}$ to all of $\Omega$. Local
super solutions can be constructed analogously. Either way we conclude the
following property of the Perron solutions $u^P_t$ corresponding to the class
$\mathcal{S}_{\bar{u}_t}$ (but cf. Remark \ref{rem: how something can be avoided}).

\begin{lemma} \label{lem: Plateau Perron stays away from boundary} For every $x
\in \partial_1 \Omega$ there exists $r(x) > 0$ such that for all sufficiently
small $t>0$, $\big(\Omega \cap \B(x, r) \big) \times (- \infty ,
\frac{\chi}{2 t})$ lies below the graph of the Perron solution $u_t^P$
corresponding to the class $\mathcal{S}_{\bar{u}_t}$, where $\bar{u}_t$ is as
in Lemma \ref{lem: barriers for Plateau problem}. Similarly, for every $x \in
\partial_2 \Omega$ there exists $r = r(x) > 0$ such that the half-infinite
cylinder $\big(\Omega \cap \B(x, r) \big) \times (- \frac{\chi}{2 t},
\infty)$ lies above this graph for all $t>0$ sufficiently small. \end{lemma}

\noindent Note that the Perron sub and super solutions described in Lemma
\ref{lem: barriers for Plateau problem} match up along $B_t:=\{ (x,
\frac{C}{\delta t^2} \dist_\Gamma(x)) : x \in \partial \Omega \text { and } -
\frac{\chi \delta t}{2 C} \leq \dist_\Gamma(x) \leq \frac{\chi \delta t}{ 2
C} \}$. It follows that the functions $u_t^P$ extend continuously from $\Omega$
to $\Omega \cup \{x \in \partial \Omega: - \frac{\chi \delta t}{2 C} \leq
\dist_\Gamma(x) \leq \frac{\chi \delta t}{ 2 C} \}$. Let $G_t := \{(x,
u_t^P(x)) :  x \in \Omega \}$ denote the graph of the Perron solution $u_t^P$
and observe that its mean curvature is pointwise bounded by $t |u^P_t| +
|\tr(p)(u_t^P)| \leq 2C$. It follows that $G_t \in \mathcal{F}_{2C}$ on the
open cylinder $\Omega \times \mathbb{R}$. From Example \ref{ex: graphs of
bounded mean curvature} in Appendix \ref{sec: almost minimizing property}
we obtain bounds for the $\mathcal{H}^{n}$-measure of $G_t$. It follows
readily that the $G_t$ are integer rectifiable $n$-varifolds in $M^n
\times \mathbb{R}$. We summarize their basic distributional properties in
the next lemma.

\begin{lemma} Let $u_t^P$ be the Perron solution for $\L_t$ on $\Omega$
constructed from the class $\mathcal{S}_{\bar{u}_t}$ where $\bar{u}_t$ is as in
Lemma \ref{lem: barriers for Plateau problem}, and let $G_t := \{(x, u_t^P(x))
:  x \in \Omega \}$ denote its graph. Then $G_t$ is an n-dimensional integer
multiplicity varifold in $M^n \times \mathbb{R}$.  In the open cylinder $\Omega
\times \mathbb{R}$, $G_t$ is a $\mathcal{C}^{2, \mu}_{\text{loc}}$-graph with
mean curvature uniformly bounded by $2 C$. We have locally uniform mass bounds
$\sup_{t \in (0, 1)} ||G_t||(W) \leq C(W) < \infty$ for every precompact set $W
\subset M^n \times (- \frac{\chi}{2t}, \frac{\chi}{2t})$. Denoting $B_t
:= \{ (x, \frac{C}{\delta t^2} \dist_\Gamma(x)) : x \in \partial \Omega \text {
and } - \frac{\chi \delta t}{2 C} \leq \dist_\Gamma(x) \leq \frac{\chi
\delta t}{ 2 C} \}$ we have that $||G_t||(B_t) = 0$, $B_t \subset \spt(G_t)$, and $\spt(G_t) \setminus
B_t$ is relatively closed in $(M^n \times (- \frac{\chi}{2t},
\frac{\chi}{2t})) \setminus B_t$.
\end{lemma}

\noindent The regularity of $G_t$ up to and including $B_t$ is immediate from Allard's boundary regularity theorem \cite{AllardBdry} once we establish that, as a varifold, $G_t$ has density $\frac{1}{2}$ along $B_t$. This is an elementary consequence of the fact that $G_t$ is a graph in $\Omega \times \mathbb{R}$, and that $B_t$ is never vertical for $t>0$. Since we are interested in the limit as $t \searrow 0$ we will supply a second, more portable argument to get boundary density $\frac{1}{2}$. We mention the abstract boundary regularity theory for almost minimizing currents in \cite{DuzaarSteffenBoundary} where some of the standard arguments below are discussed in great generality. In our concrete context we can argue directly following Lemma 5.2 in the paper \cite{AllardBdry}, whose notation we use and which all references in the proof of the following lemma are to.
For ease of notation, we will assume that $(M^n, g) = (\mathbb{R}^n, \delta)$. The author is grateful to Theodora Bourni for helping with some of the statements in \cite{AllardBdry}. Her recent work \cite{Bourni} is relevant for relaxing the regularity of the boundary $\Gamma^{n-2}$ needed in the proof.  

\begin{lemma} [Smoothness up to the boundary] \label{lem: smoothness of G_t up
to the boundary} $G_t \cap (M^n \times (- \frac{\chi}{2t},
\frac{\chi}{2t}))$ is a smooth manifold with boundary $B_t$. For every
precompact open set $W \subset M^n \times (- \frac{\chi}{2t},
\frac{\chi}{2t})$ we have that $\sup_{t \in (0, 1)} ||G_t||(W) + ||\delta
G_t||(W) \leq C(W) < \infty$.
\begin{proof}
Note that the first variation measure $\delta G_t$ of $G_t$ is a Radon measure on the complement of $B_t$ that is absolutely continuous with respect to $||G_t|| = \mathcal{H}^{n} \lfloor G_t$ with density bounded by $2C$. Since $||G_t||(B_t) = 0$, it follows from Allard's Lemma 3.1 that $\delta G_t$ extends to a Radon measure on $\mathbb{R}^n \times (- \frac{\chi}{2t}, \frac{\chi}{2t})$ such that $(\delta G_t)_{\sing} = \vec{\nu}_t d ||\delta G_t||$, where $\vec{\nu}_t$ is a $||\delta G_t||$-measurable $\mathbb{R}^{n+1}$-valued function supported on $B_t$ and normal to $B_t$ at $||\delta G_t||$ almost every point of $B_t$. By Allard's Theorem 3.5 the density function $\Theta^n (||G_t||, \cdot)$ is everywhere defined and is bounded below by $\frac{1}{2}$ at points in $B_t$. By Allard's integral compactness theorem the collection of varifold tangents at a boundary point $b \in B_t$ is non-empty and consists of integer multiplicity varifolds $C$ supported in the closed halfspace associated to
$\Tan(\partial (\Omega \times \mathbb{R}), b)$ to the side of $\Omega \times \mathbb{R}$ with the following properties:
\begin{enumerate}
\item $\frac{1}{2} \leq \Theta^n(G_t, b)= \Theta^n(||C||, 0) = \frac{||C||
(\B(0, r))}{\omega_n r^n} $ for every $r > 0$
\item $||C||(W^{\perp}) = 0$
\item $\Theta^{n}(||C||, x) \geq 1 \text{ for } ||C||-\text{a.e. } x \in
\mathbb{R}^{n+1} \setminus W^{\perp}$
\item $||\delta C|| (\mathbb{R}^{n+1} \setminus W^{\perp}) = 0 \text{ and }
\delta C \text{ extends to a Radon measure on } \mathbb{R}^{n+1}$
\item $\delta C = \vec{\nu}_C d ||\delta C|| \text{ where } \vec{\nu}_C \text{
is } ||\delta C||-\text{a.e. perpendicular to } W^{\perp}$
\end{enumerate}
Here, $W$ denotes the ($2$-dimensional) orthonormal complement of $\Tan(B_t,
b)$ in $\mathbb{R}^{n+1}$. Combining these properties with a standard monotonicity argument for varifolds of bounded first variation (as in \S 17 of \cite{LeonGMT}), it follows that any such $C$ is invariant under homotheties and hence is indeed a tangent cone. In fact, by Allard's Lemma 5.1, the Constancy Theorem for stationary varifolds, and the fact that the support of $C$ lies in a halfspace, it follows that any such tangent cone $C$ must be a finite union of integer multiplicity halfplanes through $W^{\perp}$. Allard's Theorem implies boundary regularity if $C$ consists of a single multiplicity one halfplane, or equivalently, if $\Theta^n(||G_t||, b) = \frac{1}{2}$.  Note that by Lemma \ref{lem: in F_C, current convergence is varifold convergence} and the fact that $C$ is a relative boundary in $\Omega$, the sheets composing $C$ all have multiplicity one. Away from $W^{\perp}$ the halfplanes constituting $C$ are a smooth limit of graphs of bounded mean curvature (by Allard's interior regularity theorem \cite{AllardReg}). Since $W^{\perp}$ is not vertical, this rules out the possibility of multiple sheets. A slightly more robust argument to argue that $\Theta^n(||G_t||, b) = \frac{1}{2}$, which also works when we analyze the boundary regularity in the limit as $t \searrow 0$, goes as follows: on compact subsets away from $W^{\perp}$, the tangent cone $C$ is a smooth limit of graphs whose mean curvature tends to zero. By the mean value theorem, it follows that the sheets of $C$ have to alternate in orientation (if we consider their trace in $W$) and that there must be an odd number of them. If there were three or more of them, then two consecutive oppositely oriented halfplanes would meet at an angle $< \pi$ in $W^{\perp}$. This would contradict the almost minimizing property of the approaching graphs (essentially by the triangle inequality and an argument as in Lemma \ref{lem: cylindrical almost minimizers}). Hence $C$ consists of a single, multiplicity one halfplane, as required. The other statement in the lemma follows easily now.  
\end{proof}
\end{lemma}

\begin{theorem} The family $\{G_t \lfloor \big(M^n \times (-
\frac{\chi}{2t}, \frac{\chi}{2t})\big) \}_{t \in (0, 1)}$ is precompact
in the class of integer multiplicity varifolds with locally uniformly bounded
mass and first variation in $M^n \times \mathbb{R}$, and also in the class of $2C$-almost minimizing boundaries in $\Omega \times \mathbb{R}$. In dimensions $2 \leq n \leq 7$ the possible limiting
varifolds $G$ as $t \searrow 0$ are smooth oriented hypersurfaces in $M^n \times \mathbb{R}$ with boundary $B := \Gamma^{n-2} \times \mathbb{R}$.  If $n \geq 8$ then possible limits  $G$ are smooth near $B$ and have closed singular set $\sing (G) \subset \Omega$ of Hausdorff dimension $\leq n-8$. Every connected component $G_0$ of the regular set of $G$ satisfies the marginally outer trapped surface equation $\mc_{G_0} +
\tr_{G_0} p = 0$ for one of the two possible choices for the sign of the mean
curvature of $G_0$, and $G_0$ is either a vertical cylinder or entirely graphical. By
passing to a subsequential limit of vertical translates of $G$
one obtains a cylinder $\Sigma^{n-1}\times\mathbb{R}$ with all of the above properties. Its cross-section $\Sigma^{n-1}$ is a solution of the Plateau problem for $\Gamma^{n-2}$ in $\Omega$.

\begin{proof} For the relevant compactness theorem see \S 42 in \cite{LeonGMT}
or \cite {AllardReg}. On $\Omega \times \mathbb{R}$ the varifolds $G_t \lfloor
\big(M^n \times (- \frac{\chi}{2t}, \frac{\chi}{2t})\big)$ satisfy the
$2C$-almost minimizing property. The interior regularity of a subsequential limit $G$ follows from the discussion in
Appendix \ref{sec: almost minimizing property} and Remark \ref{rem: better regularity for graphs}. The boundary regularity of $G$ follows from the `alternative argument' in the proof of Lemma \ref{lem: smoothness of G_t up to the boundary}. Lemma \ref{lem: Plateau Perron stays away from boundary} shows that away from $\Gamma^{n-2} \times
\mathbb{R}$ the limit $G$ stays away from $\partial \Omega$. Every connected
component $G_0$ of the regular set of $G$ is a smooth limit of graphs whose defining
functions $u_t$ satisfy $(\mc + \tr(p))u_t = tu_t$. The assertion about $G_0$ follows from this and the Harnack
principle in Lemma \ref{Harnack principle}. \end{proof}
\end{theorem}

\begin{remark} [Better regularity for limits of graphs] \label{rem: better
regularity for graphs} The interior regularity in dimension $n = 7$ claimed in
Theorem \ref{thm: Plateau problem} is a consequence of
the following well-known argument (see \cite{BombieriGiusti},
\cite{LeonRCE}, \cite{LeonBers}): the submanifolds we are considering here are the graphs $G_t$
of (smooth) solutions $u_t$ of the approximate Jang's equation $(\mc +
\tr(p))u_t = t u_t$ and (varifold) limits thereof. From Section \ref{sec:
interior gradient estimate} we know that the functions $\frac{1}{v_t} = (1 +
|\D u_t|^2)^{-\frac{1}{2}} > 0$ satisfy a geometric differential inequality of
the form $\Delta_{G_t} \frac{1}{v_t} - \frac{\beta}{v_t} \leq 0$, where $\beta$ is a
constant independent of $t>0$. Observe that $\frac{1}{v_t}$ is just the
vertical component of the upward pointing unit normal of $G_t$. Suppose now
that $n = 7$, and that a limit $G$ of $G_t$ has interior singular points. The
tangent cone analysis of Appendix A in \cite{LeonGMT} shows that we can then
obtain an oriented tangent cone $C$ of $G$ that is area minimizing and has an
isolated non-removable singular point at the origin. By the result of \S 2 in
\cite{BombieriGiusti} we have that $\reg(C) = C \setminus \{ 0 \}$ is
connected; by Allard's theorem, it is a smooth limit of rescalings of the
$G_t$. It hence makes sense to talk about the upward pointing unit normal of
$C$, and we conclude that its vertical component is non-negative and
superharmonic (the linear part scales away as we pass to the cone). Hence by
the Hopf maximum principle, it is either strictly positive or vanishes
identically. In the former case, $C \setminus \{0 \}$ is a minimal graph with
an isolated singularity at the origin which would then be removable by a result
of Finn's (see \cite{Finn} and also the remark in \cite{LeonRCE} at the end of
\S1), while in the latter case, $C$ is cylindrical and hence its singular set
is a vertical line. In both cases we arrive at a contradiction to the
assumption that $C$ has a non-removable isolated singularity at the origin.
\end{remark}

\begin{remark} \label{rem: how something can be avoided} Using a simple modification of these arguments (involving an application of the half-space lemmas in 36.5 and 37.6 of \cite{LeonGMT} to establish regularity) we could have avoided the use of Lemma \ref{lem: Plateau Perron stays away from boundary} here entirely, and concluded that $\partial \Omega \cap \spt \Sigma^{n-1} = \Gamma^{n-2}$ a posteriori from the strong maximum principle. 
\end{remark}

\appendix

\section{Review of geometric measure theory} \label{sec: almost minimizing
property}

In this appendix we review several classical and well-known results from geometric measure theory which have been used freely throughout this paper. The regularity theory for the classes $\mathcal{F}_C$ of `$C$-almost minimizing boundaries' is an immediate consequence of Allard's regularity theorem \cite{AllardReg} and dimension reduction as in \S 23 of \cite{LeonGMT}. 
For convenient reference we give an indication of the proofs here, following \cite{LeonGMT} closely. We confine the presentation to Euclidean space for simplicity. However, the results stated in this appendix carry over to general Riemannian manifolds and also localize, see Remark \ref{rem: localize and on Riemannian}. \\

\noindent The concept of $C$-almost minimizing boundaries below arises
naturally in the context of this paper and is convenient and appropriate for
our applications. Extensively studied and significantly more general notions include
Almgren's $(\varepsilon, \delta)$-minimal sets, see \cite{Almgren76}, and
Bombieri's $(\Psi, \omega, \delta)$-currents, see \cite{Bombieri}, and allow a
more general lower order term in (\ref{eq: defining property of C-almost
minimizing}) (typically of the form $r^{n + \alpha}$ where $r = \diam \spt
(X)$). The papers \cite{Tamanini}, \cite{Tamanini84} develop the regularity theory for almost minimizing boundaries from De Giorgi's method. After this work was finished the author has learnt that the particular almost minimizing notion we used here has been investigated systematically as that of ``$\lambda$-minimizing currents" in \cite{DuzaarSteffen}. 
We refer the reader to that paper for complete proofs and several more general results than those stated here, but we choose to keep the terminology `$C$-almost minimizing' here to emphasize that these boundaries do not come from a variational principle in our work.

\begin{definition} [C-almost minimizing property, \cite{DuzaarSteffen}] \label{def: C-almost
minimizing} A boundary $T = \partial E$ of a set $E \subset \mathbb{R}^{n+1}$
of locally finite perimeter is said to have the C-almost minimizing property if
for every open $W \subset \subset \mathbb{R}^{n+1}$ and every integer
multiplicity current $X \in \mathcal {D}_{n+1} (\mathbb{R}^{n+1})$ with support
in $W$ one has that \begin{equation} \label{eq: defining property of C-almost
minimizing} M_W(T) \leq M_W(T + \partial X) + C M_W(X). \end{equation} The
collection of all such boundaries is denoted by $\mathcal{F}_C$.
\end{definition}

\begin {remark} Let $T \in \mathcal{F}_C$, $W \subset \subset \mathbb{R}^{n+1}$
convex open, and let $Y \in \mathcal{D}_n(\mathbb{R}^{n+1})$ be integer
multiplicity with $\partial Y = 0$ and $\spt (Y) \subset W$. Then it follows
from the above definition and the sharp isoperimetric inequality (see
\cite{Almgren}) that \begin{equation*} M_W(T) \leq M_W(T + Y) + I_{n+1} C M_W
(Y) ^ {\frac {n+1}{n}} \end{equation*} where $I_{n+1}$ is the isoperimetric
constant of $\mathbb{R}^{n+1}$. Hence the $C$-almost minimizing property is more stringent than the notion studied in \cite{Tamanini}. 
\end {remark}

\begin{remark} [Bounded absolutely continuous mean curvature, cf. \S 2 in \cite{DuzaarSteffen}] \label{rem: mean
curvature bound for F_C} Definition \ref{def: C-almost minimizing} implies
that the first variation $\delta T$ corresponding to the varifold underlying $T
\in \mathcal{F}_C$ is absolutely continuous with respect to the total variation
measure $\mu_T$ associated with $T$, and $\delta T = \vec{\mc}_T d\mu_T$ for
some $\mathbb{R}^{n+1}$-valued $\mu_T$-measurable function $\vec{\mc}_T$ that
is pointwise bounded by the constant $C$. To see this, choose a variation
vector field $X \in C^1_c(W, \mathbb{R}^{n+1})$ where $W\subset
\mathbb{R}^{n+1}$ is open and bounded, and let $\phi : [0, 1] \times
\mathbb{R}^{n+1} \to \mathbb{R}^{n+1}$ denote the flow generated by $X$. Let $U
\subset W$ be an open set containing the support of $X$ such that $T \lfloor U$
is a regular slice. Using the homotopy formula (see \cite{LeonGMT}, page 139)
it follows that for $t$ small $\phi^t_\# (T) - \phi^0_\# (T) = \phi^t_\# (T) -
T = \partial \phi_\# ([0, t] \times T) = \partial \phi_\# ([0, t] \times T
\lfloor U)$ and a computation as in 26.23 in \cite{LeonGMT} shows that
$\limsup_{t \to 0} t^{-1} M_W (\phi_\# ([0, t] \times T \lfloor U)) \leq M_W(T)
\sup |X|$. Together with the $C$-almost minimizing property, the claim follows.
\end {remark}

The standard calibration argument for graphs of bounded mean curvature in the following
example was used in Section 3 of \cite{PMTII} to derive local volume bounds for graphs of
approximate solutions to Jang's equation by comparison with extrinsic balls.

\begin{example} [Graphs of bounded mean curvature] \label{ex: graphs of bounded
mean curvature} Let $u : \mathbb{R}^n \to \mathbb{R}$ be a
$\mathcal{C}^2$-function with mean curvature $|\mc_{\graph(u)} (\mathbf{x},
u(\mathbf{x}))| \leq C$ for every $\mathbf {x} = (x_1, \ldots, x_n) \in
\mathbb{R}^n$. Then $T := \partial \{(\mathbf{x}, x_{n+1}) : x_{n+1} \geq
u(\mathbf{x}) \}$ $(= \graph (u))$ has the $C$-almost minimizing property. To
see this, let $W \subset \subset \mathbb{R}^{n+1}$ open and $X \in \mathcal
{D}_{n+1} (\mathbb{R}^{n+1})$ with $\spt(X) \subset W$ as in Definition
\ref{def: C-almost minimizing}.  Denote by $\vec \nu := (1 + |\D
u|^2)^{-\frac{1}{2}} \big(\D u, -1 \big)$ the downward pointing unit normal of
$\graph (u)$ at $(\mathbf{x}, u(\mathbf{x}))$ thought of as a vector field on
all of $\mathbb{R}^{n+1}$ (parallel in the vertical direction). Observe that
$\big(\div_{\mathbb{R}^{n+1}} \vec \nu \big) (\mathbf {x}, x_{n+1}) =
\mc_{\graph(u)} (\mathbf{x}, u(\mathbf{x}))$. Let $\sigma = dx_1 \wedge \ldots
\wedge dx_{n+1} \lfloor \vec \nu$ denote the area form of $T$ and note that
$(d\sigma)(\mathbf{x}, x_{n+1}) = \mc_{\graph(u)}(\mathbf{x}, u(\mathbf{x}))
dx_1 \wedge \ldots \wedge dx_{n+1}$. Then
\begin{eqnarray*} M_W (T + \partial X) = \sup_{\omega \in \mathcal{D}^{n}(W),
|w| \leq 1} (T + \partial X) \omega  
                                    \geq  \sup_{\phi \in \mathcal{D}(W, X),
|\phi| \leq 1} (T  + \partial X) (\phi \sigma) \\
                                     \geq  \sup_{\phi \in \mathcal{D}(W, X),
|\phi| \leq 1} T (\phi \sigma) - \sup_{\phi \in \mathcal{D}(W, X), |\phi|
\leq 1} \partial X (\phi \sigma) 
                                     =     M_W (T) - \big| X (d\sigma) \big| \\ \geq  M_W(T) - C M_W(X)
\end{eqnarray*}
In this computation, we use $\mathcal{D}(W, X)$ to denote the collection of all
smooth functions $\phi$ on $\mathbb{R}^{n+1}$ with support in $W$ which are
identically equal to one in a neighborhood of the support of $X$. Hence
$T \in \mathcal{F}_C$, as asserted.
\end{example}

The following elementary consequence of Definition \ref{def: C-almost minimizing} is needed in the paper: 

\begin{lemma} \label{lem: cylindrical almost minimizers} Suppose that $T = \partial E$ is $C$-almost minimizing in $\mathbb{R}^{n+1}$, and suppose that $E = E_0 \times \mathbb{R}$. Then $T_0 = \partial E_0$ is $C$-almost minimizing in $\mathbb{R}^n$. 

\begin{proof} Let $W_0 \subset \subset \mathbb{R}^n$ be open and let $X \in \mathcal {D}_{n} (\mathbb{R}^{n})$ have support in $W_0$. Let $L\gg 1$ and consider the open set $W_L = W_0 \times (- L-1, L+1)$ and the current $X_L = X \times [-L, L]  \in \mathcal {D}_{n+1} (\mathbb{R}^{n+1})$, which has support in $W_L$. Note that $\partial X_L = \partial X \times [-L, L] + X \times \{ L\} - X \times \{ - L\}$, that $M_{W_L} ( X_L ) = 2L M_{W_0} (X)$, $M_{W_L} ( T ) = 2(L+1) M_{W_0} (T_0)$ and that $M_{W_L} (T + \partial X_L) \geq 2L M_{W_0} (T_0 + \partial X)$. The claim follows when dividing by $L $ and passing $L \to \infty$.   
\end{proof}
\end{lemma}

\begin{lemma} \label{lem: in F_C, current convergence is varifold convergence}
Let $T_j = \partial E_j \in \mathcal{F}_C$ be such that $T_j \rightharpoonup T$
as currents. Then $T = \partial E \in \mathcal{F}_C$ where $\chi_{E_j} \to
\chi_E$ in $BV_{\text{loc}}$. In fact the underlying varifolds converge as Radon
measures on $G_n(\mathbb{R}^{n+1})$.
\begin{proof} The proof is by a minor modification of the proofs of Theorems 34.5 and 37.2 in \cite{LeonGMT}, cf. also Proposition 5.1 in \cite{DuzaarSteffen}. We reproduce the argument with the necessary changes for clarity and ease of reference: it is clear from Definition \ref{def: C-almost minimizing} that $\sup_{j \geq 1} M_W (T_j) < \infty$ for every bounded open set $W$. It follows that the $\chi_{E_j}$
subconverge as $BV_{\text{loc}}$ functions to some $\chi_E$. We continue to
denote the converging subsequence by $\{j\}$. In particular, $\chi_{E_j} \to
\chi_{E}$ in $\mathcal{L}^1_{\text{loc}}$ and hence $T = \partial
E$ as currents. Next let $K \subset \mathbb{R}^{n+1}$ be compact, and let $\phi:
\mathbb{R}^{n+1} \to [0, 1]$ be a smooth function
identically equal to one near $K$ and vanishing outside an
$\varepsilon$-neighborhood of $K$. For $\lambda \in [0, 1)$ denote by $W_\lambda :=
\{ x : \phi(x) > \lambda \}$ the corresponding super level set of $\phi$. Let $R_j := E - E_j$ and note that $M_{W_0} (R_j)
\to 0$ as $j \to \infty$. Slicing these currents with respect to $\phi$ (see \S 28 of \cite{LeonGMT}), for almost every $\alpha
\in (0, 1)$ one can pass to a further subsequence (which we continue to denote
by $\{ j\}$) so that \begin{eqnarray*} && \partial (R_j \lfloor W_\alpha) =
(\partial R_j) \lfloor W_\alpha + P_j \\ && M_{W_0} (T_j \lfloor \partial
W_\alpha) = 0  \text{ for all } j \text{ and } M_{W_0}(T\lfloor \partial
W_\alpha) = 0. \end{eqnarray*} where $P_j \in \mathcal{D}_n(\mathbb{R}^{n+1})$
are integer multiplicity with support in $\partial W_\alpha$ and $M_{W_0} (P_j)
\to 0$. Hence $T \lfloor W_\alpha = T_j \lfloor W_\alpha + \partial
\tilde{R}_j + P_j$ where $\tilde{R}_j := R_j \lfloor W_\alpha$ and $P_j$ are
integer multiplicity with support in $\bar{W}_\alpha$ and whose total mass
tends to zero. For a given integer multiplicity current $X \in \mathcal{D}_{n+1} (\mathbb{R}^{n+1})$ with $\spt(X) \subset K$ the $C$-almost minimizing property
of the $T_j$'s implies that for every $\lambda \in (0, \alpha)$ one has \begin{eqnarray}
\label {eq: ugly estimate} \nonumber M_{W_\lambda} (T+\partial X) &=&
M_{W_\lambda} (T_j + \partial X + \partial \tilde{R}_j + P_j) \\ &\geq&
M_{W_\lambda} (T_j) - M_{W_\lambda}(P_j)- C (M_{W_\lambda}(X) +
M_{W_\lambda}(\tilde{R}_j)). \end{eqnarray} Letting $j \to \infty$ in (\ref{eq:
ugly estimate}) and using the lower semi-continuity of mass under weak
convergence one has \begin{equation*} M_{W_\lambda}(T + \partial X) \geq
M_{W_\lambda} (T) - C M_{W_\lambda} (X) \end{equation*} and hence $T \in \mathcal{F}_C$ (since $K$ was arbitrary). From inequality (\ref{eq: ugly estimate}) with
$X = 0$ it follows that \begin{equation*} \limsup_j \mu_{T_j} (K) \leq \limsup_j
M_{W_\lambda} (T_j)\leq M_{\{x : d(K, x) < \varepsilon \}} (T). \end{equation*}
Upon taking $\varepsilon$ to zero (more subsequences!) this gives \begin
{equation*} \limsup_j \mu_{T_j}(K) \leq \mu_T(K) \end{equation*} and hence $\mu_{T_j} \to \mu_T$ as Radon measures on $\mathbb{R}^{n+1}$. As in \cite{LeonGMT}, this argument can be repeated verbatim to show that every subsequence of the original sequence $\{T_j\}$ has a further subsequence which converges to $T$ as in the statement of the lemma. The  convergence of the original sequence follows. The first variation of the varifolds underlying the currents $\{T_j\}$ is uniformly bounded, see Remark \ref{rem: mean curvature bound for F_C}. By Allard's compactness theorem for integer rectifiable varifolds, see \S 42 of \cite{LeonGMT}, these varifolds subconverge as Radon measures on the full Grassmann
manifold $G_n(\mathbb{R}^{n+1})$ to an integer rectifiable varifold. We have already seen that their total variation measures converge to that of $\partial E$ on $\mathbb{R}^{n+1}$. Hence the varifold limit
equals the boundary $\partial E$ with its usual (multiplicity one almost everywhere) structure.
\end{proof}
\end{lemma}

\begin{lemma} [cf. Proposition 5.1 in \cite{DuzaarSteffen}] The class $\mathcal{F}_C$ is closed under translation and dilation by $\lambda \geq 1$. If $T \in \mathcal{F}_C$ then $\Theta^n(\mu_T, x) = 1$ for $\mu_T$-almost every $x \in \spt(T)$. $\mathcal{F}_C$ is sequentially compact with respect to current convergence which - for this class - is the same as varifold convergence. The tangent varifolds of $T$ are stationary and homogeneous of degree zero (and hence are tangent cones); viewed as currents, they are area minimizing. 

\begin{proof} Note that if  $T \in \mathcal{F}_C$ then $\lambda (T - x) \in \mathcal{F}_{C/\lambda}$ for every $x \in \mathbb{R}^{n+1}$ and $\lambda > 0$. The result is immediate from Remark \ref{rem: mean curvature bound for F_C}, Lemma \ref{lem: in F_C, current convergence is varifold convergence}, and a standard application of the monotonicity identity as in \S 17 of \cite{LeonGMT}.  \end{proof}  \end{lemma}

The definition of the singular and regular sets and the study of their size by dimension reduction, tangent cone analysis, and Allard's regularity theorem is exactly as in \S 36, \S 37 of \cite{LeonGMT}: 
\begin{definition} For $T \in \mathcal{F}_C$ define the regular set $ \reg (T)
:= \{ x \in \spt (T) : \exists \rho > 0 \text { with } T \lfloor \B(x, \rho)
\text{ is a connected } \mathcal{C}^{1, \alpha} \text {- graph} \}$,  and the
singular set as its complement: $\sing (T) := \spt (T) - \reg (T)$.
\end{definition}

\begin{theorem} [Regularity, cf. Theorem 1 in \cite{Tamanini} and Theorem 5.6 in \cite{DuzaarSteffen}] \label{thm: regularity of C-almost minimizers} Let
$T \in \mathcal{F}_C$. Then $\dim ( \sing (T)) \leq n - 7$. If $n=7$ the
singular set consists of isolated points.
\begin{proof} By Allard's theorem and the uniform bound on the mean curvature
for members of $\mathcal{F}_C$, see Remark \ref{rem: mean curvature bound for
F_C}, there exists $\delta > 0$ such that $\sing (T) = \{ x \in \spt (T) :
\Theta^n(\mu_T, x) \geq 1 + \delta \}$. It is not difficult to see that Leon
Simon's abstract dimension reduction argument (see Appendix A in
\cite{LeonGMT}) applies in the present context: one only has to verify that if
$T_i \rightharpoonup T$ and $x_i \in \spt (T_i)$ converges to $x \in \spt(T)$
then $\Theta^n (T, x) \geq \limsup_{i \to \infty} \Theta^n(T_i, x_i)$. This
follows easily because no mass is lost under weak convergence (i.e. $\mu_{T_i}
\to \mu_T$) and because there is an approximate monotonicity formula (see 17.4
in \cite{LeonGMT}) for the members of $\mathcal{F}_C$ because of the uniform
bound $|\vec \mc_T| \leq C$. Note that since $\mathcal{F}_C$ consists
of boundaries, the multiplicity on the regular set is always equal to $1$.
\end{proof}
\end{theorem}

\begin{remark} \label{rem: localize and on Riemannian} The
notion of being $C$-almost minimizing localizes and Theorem \ref{thm:
regularity of C-almost minimizers} continues to hold: if $\mathcal{O} \subset
\mathbb{R}^{n+1}$ is an open set, then the arguments sketched here
carry over verbatim to the class of boundaries $T = \partial E$ (in
$\mathcal{O}$) for which (\ref {eq: defining property of C-almost minimizing})
holds for all currents $X$ as described in Definition \ref{def: C-almost
minimizing} but whose support also lies in $\mathcal{O}$. The results stated here and their proofs carry over to Riemannian manifolds by isometrically embedding the manifold into some high
dimensional Euclidean space, as in \S 37 of \cite{LeonGMT}.
\end{remark}

\section {Remark on the outermost marginally outer trapped surface}
\label{sec: remark on area of outermost MOTS}

In \cite{AM} the existence, regularity, and uniqueness of the outermost marginally outer trapped surface $\Sigma^2$ (i.e. the apparent horizon) of a complete asymptotically flat $3$-dimensional initial data set $(M^3, g, p)$ is proven. It is shown that $\Sigma^2$ appears as the boundary of the total outer trapped region $T \subset M^3$. J. Metzger pointed out in \cite{Jan}  that one can force Jang's equation to blow up along $\Sigma^2$.  We remark that the method of the present paper can be used directly to conclude from such a blow-up that the apparent horizon $\Sigma^2$ has the $3|p|_{\mathcal{C}(\bar \Omega)}$-almost minimizing property on the set $\{x \in M^3 : \dist(x, \Sigma^2) < \varepsilon \} \cup M^3 \backslash T$ for some $\varepsilon > 0$. It follows that $\Sigma^2$ is \emph{outer} $3|p|_{\mathcal{C}(\bar \Omega)}$-almost minimizing, i.e. the defining property (\ref{eq: defining property of C-almost minimizing}) holds for all top dimensional currents $X \in \mathcal {D}_{3} (M^3)$ with compact support contained in the (closed) complement of $T$. \\

In the construction of the outermost marginally outer trapped surface in \cite{AM} a delicate surgery procedure is used to derive an estimate on the ``outer injectivity radius" for a certain class of marginally outer trapped surfaces $\Sigma^2 \subset M^3$ and, as a consequence, area estimates for these surfaces. The a priori curvature estimates used in this surgery procedure come from the Pogorelov-type estimate in \cite{AM05}, which relies on the Gauss-Bonnet formula. The area estimate which arises from the $C$-almost minimizing property observed in the present paper implies the area estimate in Theorem 1.2 of \cite{AM} (the assumptions there can be used to bound the area of a comparison surface) and is available in all dimensions. In \cite{GAH} we explain how the low order property of marginally outer trapped surfaces in this work brings in many classical techniques from the calculus of variations and readily leads to the existence and regularity of the apparent horizon in dimensions $n \leq 7$, without recourse to the delicate surgery procedure and stability based estimates of \cite{AM}.

\end{document}